\documentclass[11pt]{amsart}
\usepackage{verbatim, latexsym, amssymb, amsmath}
\usepackage{epsfig}
\def\R{\mathbb R}
\def\N{\mathbb N}
\def\Z{\mathbb Z}
\def\C{\mathbb C}

\def\H{\mathcal H}
\def\L{\mathcal L}
\def\x{\mathbf x}
\def\z{\mathbf z}

\def\k{\mathbf k}
\def\vol{\mathrm{vol}}
\def\re{\mathrm{Re}}

\def\a{\,d\mathcal{H}^n}
\def\b{\,d\mathcal{H}^{n-1}}

\def\MCF{mean curvature flow }

\newtheorem*{thma}{Theorem A}
\newtheorem*{thmb}{Theorem B}

\newtheorem{thm}{Theorem}[section]
\newtheorem{lemm}[thm]{Lemma}
\newtheorem{cor}[thm]{Corollary}
\newtheorem{prop}[thm]{Proposition}
\theoremstyle{remark}
\newtheorem{rmk}[thm]{Remark}
\theoremstyle{definition}

\title{Singularities of Lagrangian Mean Curvature Flow: Zero-Maslov class case}

\author{Andr\'e Neves} %${}^{\dagger}$}
\email{aneves@math.stanford.edu}
\address{Instituto Superior Te\'cnico, Lisbon, Portugal, and }
\address{Mathematics Department, Stanford University,
Stanford, CA 94305, USA}

\begin{document}

\begin{abstract}
We study singularities of Lagrangian mean curvature flow in $\C^n$ when the initial condition is a zero-Maslov
class Lagrangian. We start by showing that, in this setting, singularities are unavoidable. More precisely, we
construct Lagrangians with arbitrarily small Lagrangian angle and Lagrangians which are Hamiltonian isotopic to
a plane that, nevertheless, develop finite time singularities under mean curvature flow.

We then prove two theorems regarding the tangent flow at a singularity when the initial condition is a
zero-Maslov class Lagrangian. The first one (Theorem A) states that that the rescaled flow at a singularity
converges weakly to a finite union of area-minimizing Lagrangian cones. The second theorem (Theorem B) states
that, under the additional assumptions that the initial condition is an almost-calibrated and rational
Lagrangian, {\em connected} components of the rescaled flow converges to a {\em single} area-minimizing
Lagrangian cone, as opposed to  a possible non-area-minimizing union of area-minimizing Lagrangian cones. The
latter condition is dense for Lagrangians with finitely generated $H_1(L,\Z)$.

\end{abstract}
\maketitle \markboth{Singularities of Lagrangian Mean Curvature Flow: Zero-Maslov class case} {Andr\'e Neves}

\section{Introduction}

In the last few years, mean curvature flow of higher codimension submanifolds has attracted some attention. Most
of the work done has focused on finding initial conditions that assure the flow will exist for all time. For
instance, under some natural convexity assumptions on the image of the Gauss map, long time existence and
convergence results have been proved by J. Chen, J. Li, and Tian \cite{CLT}, Smoczyk \cite{smo1, smo2}, Smoczyk
and M.-T. Wang \cite{SM}, M.-P. Tsui and M.-T. Wang \cite{Wa4}, and M.-T. Wang \cite{Wa1, Wa0, Wa2}. On the
other hand, finite time singularities for mean curvature flow in the higher codimension case are not so well
understood and, reasoning in analogy with minimal surfaces, they are expected to exhibit a far more complicated
behavior than in the codimension one case.

%In the higher codimension case, mean curvature flow is expected to exhibit a far more , partly because such
%phenomenon is observed in minimal surfaces.

There is, therefore, interest in identifying initial conditions for the flow that are broad enough to admit
singularities, but restrictive enough so that the singularities are, so to speak, ``well-behaved''. A natural
candidate for such an initial condition is  {\em Lagrangian} because, when the ambient manifold is
K\"ahler-Einstein, the Lagrangian condition is preserved by mean curvature flow (see \cite{smo0}). Mu-Tao Wang
observed in \cite{Wa1} that, when the ambient manifold is Calabi-Yau, {\em almost-calibrated} Lagrangians (see
next section for the definition) cannot develop type I singularities, i.e., no sequence of rescaled flows at a
singularity can converge strongly to a homothetically shrinking solution. Later, Jingyi Chen and Jiayu Li
\cite{CL} showed that in this setting the sequence of rescaled flows converges weakly to an integer rectifiable
stationary Lagrangian varifold which is also a cone.

In this paper we study finite time singularities for zero-Maslov class Lagrangians in $\C^n$, a more general
condition than being almost-calibrated. The first result, Theorem A, states that the tangent flow at a
singularity can be decomposed into a finite union of area-minimizing Lagrangian cones. Theorem B is a more
interesting result because, assuming the initial condition is an almost-calibrated and {\em rational}
Lagrangian, it states that the Lagrangian angle converges to a {\em single} constant  on each {\em connected}
component of the rescaled flow. In particular, this implies that connected components of the rescaled flow
converge weakly to a single area-minimizing Lagrangian cone, instead of a possible non-area-minimizing union of
area-minimizing Lagrangian cones. Heuristically speaking, such property qualifies the formation of singularities
as being, so to speak, ``well behaved''. Without such behavior, it would be hopeless to expect Lagrangian mean
curvature flow to be more tractable than general higher codimension mean curvature flow. We remark that any
Lagrangian $M$ with $H_1(M,\Z)$ finitely generated can always be perturbed in order to become rational.

Assuming some rotational symmetry, we also construct zero-Maslov class exact Lagrangians that develop finite
time singularities under Lagrangian mean curvature flow. These examples include Lagrangians with arbitrarily
small oscillation of the Lagrangian angle and Lagrangians which are Hamiltonian isotopic to a plane.

The paper is organized as follows. In Section \ref{prelim} we recall some standard definitions and results that
will be useful throughout the rest of the paper. The main two results are discussed and stated in Section
\ref{state}. Examples of finite time singularities for Lagrangian mean curvature flow are given in Section
\ref{examples}. The first result, Theorem A, is proven in Section \ref{A}. In Section \ref{basic} we derive
evolution equations of some geometric quantities that will be needed in Section \ref{B}. In this section we
prove Theorem B.

The author would like to express his gratitude to Richard Schoen for all of his guidance and insight. He would
also like to thank Leon Simon and Brian White for enlightening discussions and constant availability.

\section{Preliminaries}\label{prelim}
Let $J$ and $\omega$ denote, respectively, the standard complex structure on $\C^n$ and the standard symplectic
form on $\C^n$. We consider also the closed complex-valued $n$-form given by
$$\Omega\equiv dz_1\wedge\ldots\wedge dz_n$$ and the Liouville form given by
$$
\lambda\equiv\sum_{i=1}^{n}x_idy_i-y_idx_i, \quad d\lambda=2\omega,
$$
where $z_j=x_j+iy_j$ are complex coordinates of $\C^n$.

A smooth $n$-dimensional submanifold $L$ in $\C^n$ is said to be {\em Lagrangian} if $\omega_L=0$ and this
implies that (see \cite{HL})
$$\Omega_L=e^{i\theta}\vol_L,$$
where $\vol_L$ denotes the volume form of $L$ and $\theta$ is some multivalued function called the {\em
Lagrangian angle}. When the Lagrangian angle is a single valued function the Lagrangian is called {\em
zero-Maslov class} and if
$$\cos \theta\geq \varepsilon$$ for some positive $\varepsilon$, then $L$ is said to be {\em almost-calibrated}.
Furthermore, if $\theta\equiv\theta_0$, then $L$ is calibrated by
$$\re\,\left( e^{-i\theta_0}\Omega\right)$$ and hence area-minimizing. In this case, $L$ is referred as being
{\em Special Lagrangian}.

 Likewise, we define an integral $n$-varifold $L_1$ and an integral $n$-current $L_2$
 to be {\em Lagrangian} if
$$
 \int_{L_1}\phi|\omega\wedge\eta|\a=0\quad\mbox{for all $n-2$ form }\eta\mbox{ and all smooth }\phi\in
  C^{\infty}_C(\C^n)
$$
and
$$
\int_{L_2}\phi\,\omega\wedge\eta\a=0\quad\mbox{for all $n-2$ form }\eta\mbox{ and all }\phi \in
C^{\infty}_C(\C^n)
$$
respectively. The concept of being Special Lagrangian can be easily extended to the case when $L$ is an integral
current.

 For a smooth Lagrangian, the
relation between the Lagrangian angle and the mean curvature is given by the following remarkable property (see
for instance \cite{schoen})
$$H=J\nabla \theta.$$

Let $L_0$ be a smooth Lagrangian in $\C^n$ such that, for some constant $C_0$, we have
$$
\H^n\bigl(L_0\cap B_R(0)\bigr)\leq C_0 R^n
$$
for all $R$ sufficiently large and assume that we have a solution $(L_t)_{0\leq t < T}$ to mean curvature flow
for which the second fundamental form of $L_t$ is bounded by a time dependent constant. The same argument used
in \cite{smo0} and the maximum principle for noncompact manifolds proved by Ecker and Huisken in
\cite{Ecker-Huisken} imply that the Lagrangian condition is preserved. In this case, we say that we  have a
solution to {\em Lagrangian mean curvature flow}. Moreover, if $L_0$ is also zero-Maslov class, then this
condition is preserved by the flow  and, according to \cite{smo}, the Lagrangian angles $\theta_t$ can be chosen
so that
$$
\frac{d\theta_t}{dt}=\Delta \theta_t.
$$
An immediate application of the parabolic maximum principle shows that the almost-calibrated condition is
preserved by Lagrangian mean curvature flow.

 A Lagrangian $L_0$ is said to be {\em rational} if for some real number $a$
$$\lambda\left(H_1(L_0, \Z)\right)= \{a2k\pi\,|\,k\in \Z\}.$$ Any Lagrangian having $H_1(L_0, \Z)$ finitely generated
can be perturbed in order to become rational. When $a=0$ the Lagrangian is called {\em exact}. Furthermore, if
$L_0$ is also zero-Maslov class, we will see in Section \ref{basic} that the rational condition is preserved by
Lagrangian mean curvature flow, i.e.,
$$\lambda\left(H_1(L_t, \Z)\right)= \{a2k\pi\,|\,k\in \Z\}$$
while the solution exists smoothly.

Assume now that the solution to mean curvature flow develops a singularity at the point $(x_0, T)$ in
space-time. Then
$$L^{\sigma}_s:=\sigma (L_{T+s/\sigma^2}-x_0) \quad \mbox{ for }-\sigma^2 T<s<0$$ is also a solution to Lagrangian
mean curvature flow and it is called a {\em rescaled flow}. It follows from \cite[Lemma 8]{ilmanen} that for
every sequence $(\sigma_i)$ going to infinity there is a a subsequence for which the \MCF
$$(L^{\sigma_i}_s)_{-\sigma_i^2T<s<0}$$ converges weakly to a homothetically shrinking weak solution of mean curvature flow
(Brakke flow). This solution is called {\em tangent flow} and depends on the sequence $(\sigma_i)$ taken.

%We end this section stating some definitions in geometric measure theory. We denote by $\H^k$ the
%$k$-dimensional Hausdorff measure of $\C^{n}$, by $\mathcal{D}^k(U)$ the set of all smooth compactly supported
%real-valued $k$-forms on the open set $U$, and by $G_k(U)$ the product
%$$U\times G(k,2n),$$ where $G(k,2n)$ is the set of all $k$-planes in $\C^n$.
% Given a countably
%$n$-rectifiable $\H^n$-measurable subset $X$ of $\C^n$ and a locally $\H^n$-integrable positive integer valued
%function $\alpha$, the {\em integral current} $T$ is the continuous linear functional on $\mathcal{D}^n(U)$
%given by
%$$
%\int_T \eta \a:=\int_{X}\langle \eta(x), \xi(x)\rangle\alpha(x)d\H^{n},
%$$
%where $\xi(x)$ is defined almost everywhere with respect to $\H^n$ and is the volume form of the approximate
%tangent $T_x X$. Likewise, the {\em integer rectifiable $n$-varifold} or {\em integral $n$-varifold} $V$ is the
%Radon measure on $G_n(U)$ given by
%$$
%\int_V \Psi:=\int_X \Psi(T_x X)d\H^n.
%$$

\section{Statement of results}\label{state}
Let $(L_t)_{0\leq t <T}$ be a smooth solution to Lagrangian mean curvature flow in $\C^n$ satisfying, for some
constant $C_0$, the area bounds
$$
\H^n\bigl(L_0\cap B_R(0)\bigr)\leq C_0 R^n
$$
for all $R$ sufficiently large. Furthermore, assume that the flow develops a finite time singularity at time $T$
and that $L_0$ is zero-Maslov class with bounded Lagrangian angle. We denote the Lagrangian angle of a rescaled
flow $(L^i_s)_{s<0}$ by $\theta_{i,s}$. Arguing informally, the following theorem states that a sequence of
rescaled flows at a singularity converges weakly to a finite union of integral Special Lagrangian cones.
\begin{thma}
If $L_0$ is zero-Maslov class with bounded Lagrangian angle, then for any sequence of rescaled flows
$(L^i_s)_{s<0}$ at a singularity, there exist a finite set $\{\bar\theta_1,\ldots,\bar\theta_N\}$ and integral
Special Lagrangian cones
$$L_1,\ldots,L_N$$
such that, after passing to a subsequence, we have for every smooth function $\phi$ compactly supported, every
$f$ in $C^2(\R)$, and every $s<0$
$$
\lim_{i \to \infty}\int_{L^i_s}f(\theta_{i,s})\phi\a=\sum_{j=1}^N m_j f(\bar\theta_j)\mu_j(\phi),
$$
where $\mu_j$ and $m_j$ denote the Radon measure of the support of $L_j$ and its multiplicity respectively.

Furthermore, the set $\{\bar\theta_1,\ldots, \bar\theta_N\}$ does not depend on the sequence of rescalings
chosen.
\end{thma}
\begin{rmk}
\
\begin{itemize}
\item[1)] It is possible and expected that, for instance, $$\{\bar\theta_1,\bar\theta_2,
\bar\theta_3\}=\{0,\pi,2\pi\}$$ but the supports of $L_1,\,L_2,$ and $L_3$ are all the same.

\item[2)] In case $n=2$, it is well known that the support of area-minimizing cones are planes  intersecting
transversely.
\end{itemize}
\end{rmk}

Theorem A follows from combining standard ideas from geometric measure theory with the evolution equation
$$\frac{d\theta_{i,s}^2}{dt}=\Delta\theta_{i,s}^2-2|H|^2.$$  We will show that, after using Huisken monotonicity
formula \cite{huisken}, such equation implies that for all $t<0$ and all positive $R$
\begin{equation}\label{eqnovo}
\lim_{i\to\infty}\int_{-1}^t\int_{L^i_s\cap B_R(0)}|H|^2+\lvert\x^{\bot}\rvert^2\a ds=0,
%\lim_{i\to\infty}\int_{-1}^t\int_{L^i_s\cap B_R(0)}|\nabla \theta_{i,s}|^2=0
\end{equation}
where $\x$ denotes the vector determined by the point $x$ in $\C^n$ and $\x^{\bot}$ denotes the projection of
the vector $\x$ onto the orthogonal complement of $T_x L^i_s$. Hence, for almost all $s<0$ we get that for all
positive $R$
\begin{equation*}
\lim_{i\to\infty}\int_{L^i_s\cap B_R(0)}|H|^2+\lvert\x^{\bot}\rvert^2\a=0
\end{equation*}
and this implies that, after passing to a subsequence, $L^i_s$ converges weakly to a stationary integral
varifold $L$ which is also a cone. Note that so far  $L$ could be a union of three Lagrangian half-planes
meeting at angles of $2\pi/3$ along a common boundary. We now sketch briefly why such configuration cannot occur
because the proof of Theorem A consists essentially in exploiting this argument. Suppose that each of the
half-planes have Lagrangian angles $\theta_1,\theta_2,$ and $\theta_3$. Then, for all sufficiently small
$\varepsilon$, $\{|\theta_{i,s}-\theta_1|\leq \varepsilon\}$ converges to a half-plane and so
$$\lim_{i\to\infty}\H(\{\theta_{i,s}=\theta_1+\varepsilon\}\cap B_R(0))>0$$
 This is impossible
because, using the coarea formula and H\"{o}lder's inequality, we have
\begin{multline*}
    \lim_{i\to\infty}\int_{-\infty}^{\infty}\H^{n-1}\left(\{\theta_{i,s}=u\}\cap B_{R}(0)\right)du
    =\lim_{i\to\infty} \int_{L^i_s\cap B_{R}(0)}|H|\a=0.
 %\leq \sqrt{D_0R^n}\left(\int_{M^i\cap B_{2R}(0)}|H|^2\right)^{1/2}.
\end{multline*}

Theorem A raises the following question: Given $\Sigma^i$
 a sequence of connected components of $L^i_s\cap B_R(0)$ that converges weakly to $\Sigma$, does $\Sigma$ need
to be a SLag cone? In other words, does $\theta_{i,s}$ need to converge to a constant? According to Theorem A we
only know that $\Sigma$ is a finite union of Slag cones which might have different Lagrangian angles and hence
not necessarily area-minimizing. An affirmative answer to this question is necessary if one wants to make
reasonable the possibility of developing a regularity theory for the flow.

Technically, the difficulty  comes from the fact that because the sequence of smooth manifolds $L^i_s$ are
becoming singular when $i$ goes to infinity, no Poincar\'e inequality holds with a constant independent of $i$
and therefore we cannot conclude that, on each connected component of $L^i_s$, the Lagrangian angles
$\theta_{i,s}$ converge to a constant. As a matter of fact, for the sequence of smooth surfaces
$$ L_{\varepsilon}\equiv \{(z,w) \in \C^2\,|\,zw=\varepsilon\},$$
one can easily construct bounded functions $f_{\varepsilon}$ for which the $L^2$ norm of its gradient goes to
zero when $\varepsilon$ goes to zero but nevertheless $f_{\varepsilon}$  converges to a distinct constant on
each complex plane. The question raised in the previous paragraph was addressed in \cite[Theorem 5.1]{CL} but
unfortunately this technical aspect was overlooked.
%Moreover, there is no pointwise bound on $|\nabla \theta_{i,s}|$ and hence, even if its $L^2$-norm goes to zero,
%it could be possible that the Lagrangian angle on a connected component of $L^i_s$ converges to two distinct
%values.

In order to deal with this difficulty, we require $L_0$ to satisfy two additional conditions, namely that it is
an almost-calibrated and rational Lagrangian (see Section \ref{prelim} for the definitions). We argued in
Section \ref{prelim} that these conditions are preserved by Lagrangian mean curvature flow.
\begin{thmb}
If $L_0$ is almost-calibrated  and rational, then, after passing to a subsequence of $(L^i_s)_{s<0}$, the
following property holds for all $R>0$ and almost all $s <0$. \par For any convergent subsequence (in the Radon
measure sense) $\Sigma^i$ of connected components of $B_{4R}(0)\cap L^i_s$ intersecting $B_R(0)$, there exists a
Special Lagrangian cone $L$ in $B_{2R}(0)$ with Lagrangian angle $\bar\theta$ such that
$$
\lim_{i \to \infty}\int_{\Sigma^i}f(\theta_{i,s})\phi\a=m f(\bar\theta)\mu(\phi),
$$
for every $f$ in $C(\R)$ and every smooth $\phi$ compactly supported in $B_{2R}(0)$, where $\mu$ and $m$ denote
the Radon measure of the support of $L$ and its multiplicity respectively.
\end{thmb}

Next, we give a heuristic argument explaining why the rational condition should play a role.  From the
pioneering work of Richard Hamilton both on Ricci flow and on mean curvature flow we know that it is helpful to
find quantities that are constant on self-similar solutions. For that matter, let us consider $$L_s\equiv
\sqrt{s}L_1$$ to be a solution  to Lagrangian mean curvature flow where $L_0$ is zero-Maslov class. A simple
computation reveals that for all $s>0$
\begin{align*}
H(L_s)=\x^{\bot}/(2s) & \iff 2s\nabla\theta_s=-(J\x)^{\top}\\
& \iff 2s\,d\theta_s+\lambda=0.
\end{align*}
Thus, we conclude that $L_s$ is exact and that if we denote by $\beta_s$ the primitive for the Liouville form
$\lambda$, then $\beta_s+2s\theta_s$ is constant in space for all $s$. Arguing informally, this suggests that
showing convergence of the  Lagrangian angle to a single constant  should be equivalent to showing that the
primitive for the Liouville form converges to a single constant. The advantage of doing so is that the gradient
of $\beta_s$ is a first order quantity and thus easier to control than the gradient of $\theta_s$ which is a
second order quantity.

We now sketch the main idea behind the proof of Theorem B. Assume, for the sake of simplicity, that $L_0$ is
exact which implies that for each $i$ there is a family of smooth functions $\beta_{i,s}$ defined on $L^i_s$
such that $ d\beta_{i,s}=\lambda,$ or equivalently,
$$J\nabla \beta_{i,s}(x)=-\x^{\bot}\quad\mbox{for all }x \in M_{i,s}.$$
Moreover, as it will be shown in Section \ref{basic}, the functions $\beta_{i,s}$ can be chosen so that
 $$\frac{d}{ds}(\beta_{i,s}+2s\theta_{i,s})=\Delta(\beta_{i,s}+2s\theta_{i,s}).$$
This evolution equation combined with identity \eqref{eqnovo} implies that, after passing to a subsequence,
$\beta_{i,s}+2s\theta_{i,s}$ has a limit which is independent of $s$ and so it must converge to some constant
$c_j$ on each Special Lagrangian cone $L_j$, with $j=1,\ldots,N.$ Hence, we obtain from Theorem A that
$\beta_{i,s}$  converges to $c_j-2s\bar\theta_{j}$ on each $L_j$. Moreover, we can assume without loss of
generality that the set
$$\Lambda\equiv\{\bar\beta_1-2s\bar\theta_1,\ldots,\bar\beta_N-2s\bar\theta_N\}$$ has $N$ distinct values.

Let $\Sigma_i$ be a convergent sequence of connected components of $L^i_s\cap B_R(0)$.
 Because the gradient of $\beta_{i,s}$ is {\em pointwise bounded}  and its $L^2$-norm converges to zero, we can show
 that the sequence of functions $\beta_{i,s}$ converges to a single constant when restricted to $\Sigma_i$
 (see Proposition \ref{cont}). Thus, the Lagrangian angle of $\Sigma_i$ must converge to a constant
 because otherwise two numbers in the set $\Lambda$ would be equal.

\section{Examples of Finite Time Singularities}\label{examples}

We construct examples of finite time singularities for mean curvature flow where the initial condition is a
zero-Maslov class and exact Lagrangian.

For simplicity, we restrict ourselves to  $\C^2$ but we note that the phenomena observed also occur in $\C^n$.
Given a curve $\gamma$ in the complex plane, it is easy to see that
$$
L=\{(\gamma\cos \alpha, \gamma \sin \alpha)\;|\; \alpha \in \R/2\pi\Z\}
$$
is a Lagrangian surface in $\C^2$. A choice of orientation for the curve $\gamma$ induces an orientation on $L$
and if $\gamma(s)$ denotes a parametrization of $\gamma$, then
$$
\Omega_L=\frac{\gamma}{|\gamma|}\frac{\gamma^{\prime}}{|\gamma^{\prime}|}\vol_L\quad\mbox{and}\quad\lambda_L=\langle
i\gamma,\gamma^{\prime}\rangle ds.
$$
Hence, we get that $L$ is exact and zero-Maslov class whenever $\gamma$ is diffeomorphic to a line.

If we evolve $L$ by mean curvature flow, the rotational symmetries are preserved and the corresponding
$\gamma_t$ evolve according to
\begin{equation}\label{flow}
\frac{d z}{dt}=\k-\z^{\bot}/|z|^2,
\end{equation}
where $\k$ is the curvature of $\gamma$ and $\z^{\bot}$ denotes the projection of the position vector $\z$ on
the orthogonal complement of $T_{x}\gamma$.

\begin{figure}
\centering {\epsfig{file=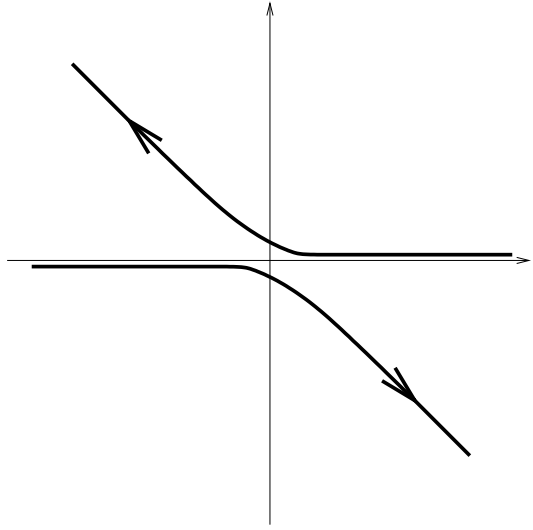, height=200pt}}\caption{Lagrangian surface $L_0$.} \label{fig2}
\end{figure}

For any $0<\beta\leq\pi$, consider the following initial condition for the {\em equivariant mean curvature flow}
\eqref{flow}
$$
\gamma_0(s)=(\sin(\pi s/\beta))^{-\beta/\pi}e^{is} \equiv r_0(s)e^{is}, \quad 0<s<\beta.
$$
The corresponding Lagrangian surface $L_0$ is asymptotic to two oriented planes with Lagrangian angles $\pi$ and
$2\beta$ and, when $\beta>\pi/2$, its intersection with $\C\times\{0\}$ can be seen in Figure \ref{fig2}.
In order to compute the Lagrangian angle of $L_0$, we use the formula
$$\theta_0(s)=\arg(\gamma_t\gamma_t^{\prime})=2s+\arg(r_0^{\prime}+ir_0)$$
and obtain that
$$\theta_0(s)=(2-\pi/\beta)s+\pi, \quad 0<s<\beta.$$
Note that the oscillation of the Lagrangian angle can be made arbitrarily small by choosing $\beta$ close to
$\pi/2$.

\begin{figure}
\centering {\epsfig{file=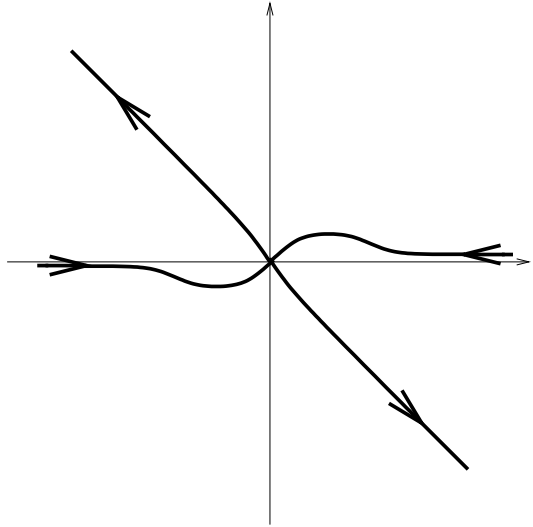, height=200pt}}\caption{Finite time singularity at the origin.} \label{fig4}
\end{figure}

We now sketch briefly three distinct behaviors for the equivariant mean curvature flow. When $0<\beta<\pi/2$,
 the curve will expand indefinitely because the curvature term on the right hand side of \eqref{flow} points outward
 and dominates the first-order term that points inward. As a matter of fact, Anciaux \cite{anciaux} found a self-expander
 with the same asymptotics at infinity as $\gamma_0$. When $\beta=\pi/2$, the Lagrangian surface is one of the
 Special Lagrangians studied in  \cite{HL}. Thus, the curvature term equals the first order term on \eqref{flow} because
  the curve is a fixed point for the flow. Finally,
  when $\pi/2<\beta\leq \pi$, the first order term will be pointing inward and bigger than the curvature term, thus
  forcing the solution to have a finite time singularity at the origin (see Figure \ref{fig4}). This is the
  content of the next theorem.

\begin{thm}\label{example}
When $\pi/2<\beta\leq \pi$, the Lagrangian mean curvature flow starting at $L_0$ develops a finite time
singularity at the origin. The tangent flow is a union of two planes intersecting at a single point, both with
Lagrangian angle $\beta/2$.
\end{thm}
\begin{proof}
We start by proving short-time existence for the equivariant mean curvature flow. The procedure is well-know
among the specialists but we include it here for the sake of completeness.

 After rotating the coordinate axis by
$(\pi-\beta)/2$, the curve $\gamma_0$ can be written as the graph of a function $u_0$ over the real axis. A
straightforward computation shows the existence of some constant $C$ such that
\begin{equation}\label{asym}
|u'_0|_{C^2}+|xu_0'-u_0|_{C^0}\leq C.
\end{equation}
For each fixed $n\in \N$, consider graphical solutions $\gamma^n_t\equiv (x,u_t^n(x))$ for the equivariant mean
curvature flow with boundary conditions
$$u^n_0(x)=u_0(x)\quad\mbox{for } |x|\leq n, \quad u^n_t(\pm n)=u_0(\pm n)= u_0(n).$$ We will show uniform apriori
$C^{2,\alpha}$-estimates for the sequence of functions $(u^n_t)$.

A simple computation reveals that $u^n_t$ solves the quasilinear equation
\begin{equation}\label{quasi}
\frac{du}{dt}=\frac{u''}{1+(u')^2}+\frac{xu'-u}{x^2+u^2}.
\end{equation}
\begin{lemm}\label{max}
There exists positive $s_0$ and $\varepsilon$ so that $$u^n_t(0)\geq\varepsilon$$ for all $t\leq s_0$ and all
$n\in\N$. Moreover, we have for all $t\leq s_0$ that
\begin{equation*}
 u_0(n)|x|/n\leq {u^n_t}(x)\leq u_0(x).
\end{equation*}
\end{lemm}
\begin{proof}
 Consider a solution $(C_t)_{t\geq 0}$ to \eqref{flow} having initial condition a circle of small radius centered
at the origin that does not intersect $\gamma_0$.  The maximum principle implies that the graph of $u^n_t$
cannot intersect $C_t$ and so the first assertion follows. The second assertion also follows from the maximum
principle because
$$v^n(x)\equiv u_0(n)|x|/n$$
and $u_0$ are a solution and supersolution for \eqref{quasi} respectively.
\end{proof}

The function $$v^n_t\equiv u^n_t-u_0$$ satisfies the equation
$$
\frac{dv}{dt}=\frac{v''}{1+(u'+u_0')^2}+\frac{xv'-v}{x^2+(v+u_0)^2}+F_t
$$
where, due to \eqref{asym},
$$F_t\equiv \frac{u_0''}{1+(u'+u_0')^2}+\frac{xu_0'-u_0}{x^2+(v+u_0)^2}$$
is pointwise bounded. Hence, the maximum principle implies that $v^n_t$ is uniformly bounded for all $t\leq
s_0$. Moreover, we obtain from Lemma \ref{max} that
$$-u_0(n)/n\leq {u^n_t}'(-n)\leq u_0'(-n)\quad\mbox{and}\quad u_0'(n)\leq {u^n_t}'(n)\leq u_0(n)/n$$
and so, it follows from \eqref{asym} that ${v^n_t}'(\pm n)$ converges to zero as $n$ goes to infinity. Because
$\phi^n_t\equiv {v^n_t}'$ satisfies an evolution equation of the form
$$\frac{d\phi}{dt}=a(x,\phi')\phi''+b(x,\phi,\phi')\phi'+c(x,\phi,\phi')\phi+G_t,$$
where $a>0$ and $c, G_t$ are uniformly bounded functions, we obtain from the maximum principle that ${v^n_t}'$
is uniformly bounded. Standard theory for quasilinear parabolic equations implies the existence of some constant
$M$ for which $|u^n_t-u_0|_{C^{2,\alpha}}<M$ for all $t\leq s_0$. Therefore, we can let $n$ go to infinity and
obtain a solution $\gamma_t(x)\equiv(x,u_t(x))$ for the equivariant mean curvature flow.

Next, we argue that the flow $(\gamma_t)$ develops a finite time singularity. We need the following lemma.
\begin{lemm}\label{moby}
While the solution exists smoothly, the curve $\gamma_t$ can be parametrized by
$$
\gamma_t(s)=r_t(s)e^{is}\quad\mbox{with }r_t(s)>0,\quad 0<s<\beta.
$$
\end{lemm}
\begin{proof}
    For any $0<\alpha<\beta$, denote by $C_{\alpha}$ the line
        $$C_\alpha=\{r\exp^{i\alpha}\,|\, r\in\R \}.$$
    Initially, we have that $C_\alpha$ and $\gamma_0$ intersect only once. Furthermore, it follows from the short-time
    existence estimates that $\gamma_t$ remains in the region below $\gamma_0$ and above the $x$-axis. Hence, the Sturmian Theorem
    proved by Angenent \cite[Proposition 1.2.]{ang} implies that $C_{\alpha}$ and $\gamma_t$ must intersect exactly once
    while the solution exists smoothly.
\end{proof}

For the rest of this proof we parameterize the curves $\gamma_t$ as described in the previous lemma. The
equation satisfied by $r_t$ becomes
\begin{lemm}\label{eqaux}
\begin{equation*}
\frac{dr}{dt}=-\frac{\theta_t^{\prime}}{r}=\frac{r r^{\prime\prime}-2r^2-3(r^{\prime})^2}{r(r^{\prime})^2+r^3},
\end{equation*}
\end{lemm}
\begin{proof}
 Denote by $\partial_s$ the tangent vector
 $$\partial_s=r'e^{is}+ire^{is}.$$
 Then,
 $$\langle d(re^{is})/ds, i\partial_s\rangle=dr/dt\langle e^{is},i\partial_s\rangle=-rdr/dt. $$
 On the other hand,
 $$\langle d(re^{is})/ds, i\partial_s\rangle=\langle H,i\partial_s\rangle=\langle \nabla\theta_t,\partial_s\rangle=
 \theta'_t$$
 and so the first identity follows. The second identity can be checked using
$$\theta_t(s)=2s+\arg(r_t^{\prime}+ir_t).$$
\end{proof}
Let $A_t(\varepsilon)$ denote the area of the triangular-shaped region
$$\{ue^{is}\,|\, \varepsilon \leq s\leq \beta-\varepsilon,\,0\leq u\leq r_t(s)\}.$$
Note that
$$2A_t(\varepsilon)=\int_{\varepsilon}^{\beta-\varepsilon}r^2_t(s)\,ds$$
and that $$2s<\theta_t(s)<2s+\pi$$ because $\theta_t(s)=2s+\arg(r_t^{\prime}+ir_t).$ Therefore,
\begin{align*}
\frac{d}{dt} A_t(\varepsilon) & =-\int_{\varepsilon}^{\beta-\varepsilon}\theta_t^{\prime}(s)\,ds =
(\theta_t(\varepsilon)-\theta_t(\beta-\varepsilon)) <\pi+2\varepsilon-2\beta.
\end{align*}
Because $\varepsilon$ can be chosen arbitrarily small, the flow must develop a finite time singularity if
$\pi/2<\beta\leq \pi$.

Denote by $T$ the instant of the first time singularity. We need to show that the singularity occurs at the
origin. The key idea consists in showing that if that is not the case, then the tangent flow cannot be a union
of Lagrangian planes, which is a contradiction to Theorem A. In order to do so, we need some preliminary lemmas.

\begin{lemm}For all $t<T$
$$\lim_{s\to 0}\theta_0(s)=0\quad\mbox{and}\quad\lim_{s\to \beta}\theta_t(s)=2\beta.$$
\end{lemm}
\begin{proof}
    The maximum
    principle applied to $\theta_t$ implies that $\pi\leq \theta_t\leq 2\beta$ for all $t<T$. Suppose that
    there is $t_1<T$, a sequence $(s_i)$
    converging to zero, and $\varepsilon>0$ so that
        $$\lim_{i\to\infty}\theta_{t_1}(s_i)=\pi+2\varepsilon.$$

    Recall that $L_t$ denotes the Lagrangian surfaces corresponding to $\gamma_t$ and consider the function
    $$\phi_{t,\varepsilon}\equiv(\theta_t-\pi-\varepsilon)^3_{+}$$
    which is supported on $\{p\in L_t\,|\,\theta_t\geq\pi+\varepsilon\}$. Furthermore,
    $$\frac{d\phi_{t,\varepsilon}}{dt}\leq\Delta\phi_{t,\varepsilon}.$$
    Huisken's monotonicity formula \cite{huisken} implies that for all $i\in\N$
    \begin{multline*}
        8\varepsilon^3\leq\int_{L_0}\phi_{t,\varepsilon}\frac{\exp(-|x-x_i|^2/4t_1)}{4\pi t_1}d\H^2\\=
        \int_{\{\theta_0\geq\pi+\varepsilon\}}\phi_{t,\varepsilon}\frac{\exp(-|x-x_i|^2/4t_1)}{4\pi t_1}d\H^2,
    \end{multline*}
    where $x_i$ is the point $(\gamma_{t_1}(s_i),0)$ in $\C^2$. For every $R>0$, we have for all $i$
    sufficiently large that
        $$\{\theta_0\geq\pi+\varepsilon\}\cap B_R(x_1)=\emptyset.$$
    Thus
        $$\lim_{i\to\infty}\int_{\{\theta_0\geq\pi+\varepsilon\}}
        \phi_{t,\varepsilon}\frac{\exp(-|x-x_i|^2/4t_1)}{4\pi t_1}d\H^2=0$$
    and this gives us a contradiction.
\end{proof}
This lemma is used to prove
\begin{lemm}\label{angle}
For all $t<T$
 $$\frac{dr}{dt}\leq0.$$
\end{lemm}
\begin{proof}
Taking into account that the parameterization  described in Lemma \ref{moby} creates a tangential component on
the deformation vector, we get that
$$
\frac{d\theta}{dt}=\Delta_{L_t}\theta+\left\langle \frac{dx}{dt},\nabla \theta_t\right\rangle
 =\frac{\theta^{\prime\prime}}{|\gamma^{\prime}|^2}+\theta^{\prime}\left(\frac{1}{r|\gamma^{\prime}|}
 \left(\frac{r}{|\gamma^{\prime}|}\right)^{\prime}+\frac{dr}{dt}\frac{r^{\prime}}{|\gamma^{\prime}|^2}
\right).
$$
While the solution exists smoothly, we have that
$$\lim_{s\to 0}\theta_t(s)=\pi\quad\mbox{and}\quad\lim_{s\to \beta}\theta_t(s)=2\beta$$ and thus,
 the Sturmian property \cite[Proposition 1.2.]{ang} implies that the cardinality
$$
\#\{s\,|\,\theta_t(s)=y\}
$$
is one if $\pi<y<2\beta$ and zero if $y<\pi$ or $y>2\beta$. Hence
$$\frac{dr}{dt}=-\frac{\theta^{\prime}_t}{r}\leq0$$
for all $t<T$.
\end{proof}
The curves $\gamma_t$ are symmetric under refection over a line with slope $\tan(\beta/2)$ and so
\begin{equation}\label{symm}
r_t(\beta/2+s)=r_t(\beta/2-s)
\end{equation}
for all $t<T$. This implies that
$$
r_t^{\prime}(\beta/2)=0\quad\mbox{for all }t<T.
$$
\begin{lemm}\label{critical}
For any  $t<T$, $r_t(s)$ is decreasing when $s<\beta/2$ and increasing when $s>\beta/2$.
\end{lemm}
\begin{proof}
Direct computation shows that $\beta/2$ is the only critical point of $r_0$ and that, denoting $r_t^{\prime}$ by
$u_t$,
$$
\frac{du_t}{dt}=\frac{u_t^{\prime\prime}}{(r^{\prime})^2+r^2}+u_t^{\prime}b(r_t,u_t,u_t^{\prime})+u_t
c(r_t,u_t,u_t^{\prime}),
$$ where the functions $b$ and $c$ are bounded for each $t<T$. Moreover,
$$\lim_{s\to 0}u_t(s)=\infty\quad\mbox{and}\quad\lim_{s\to\beta}u_t(s)=-\infty$$
and thus, the Sturmian property \cite[Proposition 1.2.]{ang}
 implies that $\beta/2$ is the only critical point of $r_t$.
\end{proof}

Suppose now that the singularity happens at a point $x_0\equiv ae^{i\alpha}$, with $0<a\leq r_0(\alpha)$ and
$0<\alpha<\beta$. From Theorem A, we know that the tangent flow at the singularity is a union of planes and so,
by White's regularity Theorem \cite{white},
$$\limsup_{\delta\to 0}\frac{\H^1\left(\gamma_{T-\delta^2}\cap B_{\delta}(x_0)\right)}{2\delta}\geq 2.$$ We
show next that this is impossible because for all $\delta$ sufficiently small and all $t<T$
$$
\frac{\H^1\left(\gamma_{t}\cap B_{\delta}(x_0)\right)}{2\delta}\leq 3/2.
$$
Without loss of generality we assume that $\alpha=\beta/2$ (the remaining cases are treated similarly). For any
$\delta<a$, Lemma \ref{angle} and Lemma \ref{critical} imply that
$$\gamma_{t}\cap B_{\delta}(x_0)$$ is either empty or  a connected curve. If the latter occurs,
there is $\varepsilon(t)<\arcsin(\delta/a)$ for which
$$
\gamma_{t}\cap B_{\delta}(x_0)=\{\gamma_t(s)\,\,|\,\,  |s-\beta/2|<\varepsilon\}.
$$
Note that
$$
(r_t(\beta/2+\varepsilon)\cos(\varepsilon)-a)^2+(r_t (\beta/2+\varepsilon)\sin(\varepsilon))^2=\delta^2
$$
and so
$$
|r_t(\beta/2+\varepsilon)\cos(\varepsilon)-a|+|r_t (\beta/2+\varepsilon)\sin(\varepsilon)|\leq \sqrt
2\delta<3/2\delta.
$$
Combining this inequality with Lemma \ref{angle}, Lemma \ref{critical}, and \eqref{symm}, we obtain
\begin{align*}
\frac{\H^1\left(\gamma_{t}\cap B_{\delta}(x_0)\right)}{2\delta} &
=\frac{1}{\delta}\int_{\beta/2}^{\beta/2+\varepsilon}((r_t^{\prime})^2+r_t^2)^{1/2}ds\\
& \leq \frac{r_t(\beta/2+\varepsilon)-r_t(\beta/2)}{\delta}+\varepsilon\frac{r_t(\beta/2+\varepsilon)}{\delta}\\
& \leq \frac{r_t(\beta/2+\varepsilon)-a}{\delta}+\varepsilon \frac{r_t(\beta/2+\varepsilon)}{\delta}\\
& \leq 3/2
\end{align*}
for all $\delta$ sufficiently small.

Finally, we argue next that the tangent flow at the singularity is a union of two planes with Lagrangian angle
$\pi/2+\beta$. From \eqref{symm} it follows that
$$
\theta^{\prime}_t(\beta/2+s)=\theta^{\prime}_t(\beta/2-s)
$$
and therefore, because the solution remains asymptotic to two planes with Lagrangian angles $\pi$ and $2\beta$,
we obtain after integration that $\theta_t(\beta/2)=\pi/2+\beta$. From Lemma \ref{critical} we know that
$\gamma_t(\beta/2)$ is the closest point of $\gamma_t$ to the origin and so Theorem B implies the desired
result.
\end{proof}

\begin{figure}
\centering {\epsfig{file=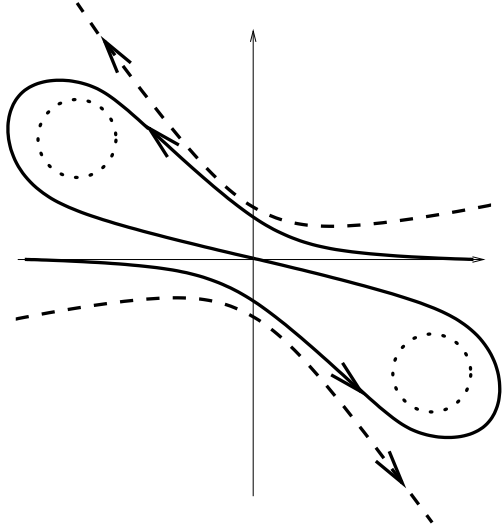, height=200pt}}\caption{Lagrangian Hamiltonian isotopic to a plane developing
a finite time singularity.} \label{fig3}
\end{figure}

We can now use Theorem \ref{example} to construct an exact and zero-Maslov Lagrangian class which is Hamiltonian
isotopic to a Lagrangian plane that, nevertheless, develops a finite time singularity. Denote by $L_0$ the
compact perturbation of a Lagrangian plane which is associated with the curve described in Figure \ref{fig3}.
The dashed noncompact curve represents one of the curves described in Theorem \ref{example} (slightly rotated so
that it is not asymptotic to $L_0$) and has a finite time singularity at the origin at time $T$. The dashed
circles shown in Figure \ref{fig3} correspond to a Lagrangian torus, which will have a finite time singularity
at time $T_1$. All these curves can be arranged so that $T<T_1$ and an explicit expression for such curves could
be easily found. The short-time existence for the flow with initial condition $L_0$ follows from the same
arguments used  in the proof of Theorem \ref{example}. Because the two noncompact solutions we consider have
different asymptotics, the maximum principle implies that they can never intersect. Hence, the flow
$(L_t)_{t\geq 0}$ must develop a finite-time singularity.

We end this section with a brief heuristic discussion of how could the flow $(L_t)_{t\geq 0}$ be continued after
its finite-time singularity.  It is expected that in the setting described above, the singularity occurs at the
origin. In this situation, the Lagrangian surface at the time of the singularity decomposes into a union of an
immersed $2$-sphere (the immersion point being at the origin) and a Lagrangian surface diffeomorphic to the
Lagrangian plane. As it was pointed out by Tom Ilmanen, there are two possible different evolutions for the
Lagrangian surface after the singularity occurs: the immersed $2$-sphere that has formed can evolve as an
immersed $2$-sphere or it can become an embedded torus which then evolves smoothly by mean curvature
 flow. In either case, the other connected piece will evolve smoothly to a Lagrangian plane.

\section{Proof of Compactness Theorem A}\label{A}
The next proposition will be essential to prove Theorem A. As a mean of motivation, it could be easier to read
first the proof of Theorem A and come back to the proposition when necessary.
\begin{prop}\label{general}
Let $(L^i)$ be a sequence of smooth zero-Maslov class Lagrangians in $\C^n$ such that, for some fixed $R>0$, the
following properties hold:
\begin{itemize}
\item[(a)] There exists a constant $D_0$ for which $$\H^n(L^i\cap B_{2R}(0))\leq
D_0R^n\quad\mbox{and}\quad\sup_{L^i\cap B_{2R}(0)}|\theta_i|\leq D_0$$ for all $i \in \N$.

 \item[(b)]
$$\lim_{i \to \infty}\H^{n-1}(\partial L^i\cap B_{2R}(0))=0$$
and $$\lim_{i \to \infty}\int_{L^i\cap B_{2R}(0)}|H|^2\a=0.$$
\end{itemize}
Then, there exist a finite set $\{\bar\theta_1,\ldots,\bar\theta_N\}$ and integral Special Lagrangians
$$L_1,\ldots,L_N$$
such that, after passing to a subsequence, we have for every smooth function $\phi$ compactly supported in
$B_R(0)$ and every $f$ in $C(\R)$
$$
\lim_{i \to \infty}\int_{L^i}f(\theta_{i})\phi\a=\sum_{j=1}^N m_j f(\bar\theta_j)\mu_j(\phi),
$$
where $\mu_j$ and $m_j$ denote, respectively, the Radon measure of the support of $L_j$ and its multiplicity.
\end{prop}

\begin{proof}
 From Allard compactness theorem for varifolds \cite[Theorem 42.7]{Leon} we obtain the existence of a subsequence, still
 denoted by $(L^i)$,
 converging
in $B_{2R}(0)$ to a stationary integer rectifiable varifold $L$. Moreover,
$$
\int_{L}\phi|\omega\wedge\eta|\a=0
$$
for every $n-2$ form $\eta$ and all smooth $\phi\in C^{\infty}_C(B_{2R}(0))$,  and this implies that $L$ is
Lagrangian. It suffices to find integral Special Lagrangians
$$L_1,\ldots,L_N,$$ a finite set $\{\bar\theta_1,\ldots,\bar\theta_N\},$ and some positive $\varepsilon_0$
such that, after passing to a subsequence of $(L^i)$, we have for all smooth $\phi$ compactly supported in
$B_R(0)$, all $0<\varepsilon<\varepsilon_0$, and all $j=1,\ldots,N$,
$$
\lim_{i \to \infty}\int_{\{|\theta_i-\bar\theta_j|\leq\varepsilon\}}\phi\a= m_j\mu_j(\phi)$$ and
$$
\mu_L(\phi)=\sum_{j=1}^N m_j\mu_j(\phi),
$$
where $\mu_L$ and $\mu_j$ denote the Radon measure of $L$ and of the support of $L_j$ respectively, and $m_j$
denotes the multiplicity of $L_j$.

The idea for the proof is as follows. The regular points of $L$ form a dense open set and therefore we can pick
$p$ in $L\cap B_R(0)$ such that, for some positive $\rho$, $B_{\rho}(p)$ is contained in $B_R(0)$ and the
support of $L\cap B_{\rho}(p)$ is a smooth Special Lagrangian with angle $\bar\theta_1$. After adding some
multiple of $\pi$ to $\bar\theta_1$ if necessary, we will show the existence of an integral Special Lagrangian
$L_1$ and of
 $\varepsilon_1>0$ such that, for all smooth $\phi$ with compact support in $B_{R}(0)$ and all
$0<\varepsilon\leq\varepsilon_1$, we have
\begin{equation}\label{final}
    \lim_{i \to \infty}\int_{\{|\theta_i-\bar\theta_1|\leq\varepsilon\}}\phi\a= m_1\mu_1(\phi),
\end{equation}
where $\mu_1$ is the Radon measure of the support of $L_1$ and $m_1$ its multiplicity. Because the support of
$L_1$ is stationary, the monotonicity formula implies that
\begin{equation}\label{tricky}
  \mu_1(B_{2R}(0))R^{-n}\geq \mu_1(B_R(p))R^{-n}\geq \mu_{1}(B_{\rho}(p))\rho^{-n}\geq \gamma_n
\end{equation}
 for some universal constant $\gamma_n$.

In order to find $\bar\theta_2$ and the integral Special Lagrangian $L_2$, we repeat this process but this time
applied to the sequence $$P_i\equiv\{|\theta_i-\bar\theta_1|\geq \varepsilon_1\},$$ where the boundary will
cause no difficulty because, as it will be seen in the proof of Lemma \ref{lemma0}, we can assume that
$$
\lim_{i\to\infty}\H^{n-1}(\{\theta_i=\bar\theta_1\pm\varepsilon_1\}\cap B_{2R}(0))=0
$$
and hence,
\begin{multline*}
    \lim_{i\to\infty}\H^{n-1}(\partial P_i\cap B_{2R}(0))\\\leq \lim_{i\to\infty}\left(\H^{n-1}(\partial L_i\cap
    B_{2R}(0))+\H^{n-1}(\{\theta_i=\bar\theta_1\pm\varepsilon_1\}\cap B_{2R}(0))\right)=0.
\end{multline*}

 Condition (a) and \eqref{tricky} ensures that this will be done only finitely many times and hence the proposition will be
proven as soon as we show \eqref{final}.

The next lemma will be quite useful throughout the rest of the proof.
\begin{lemm}\label{lemma0}
For almost all endpoints $a$ and $b$, the sequence $$N^i\equiv \{a\leq\theta_i\leq b\}$$ contains a subsequence
converging, in $B_{2R}(0)$,  to a stationary integer rectifiable varifold $N$ in the varifold sense and to an
integral current $ \widehat{N}$ with $\partial \widehat{N}=0$ in the current sense.
\end{lemm}

\begin{proof}
For almost all endpoints $a$ and $b$ we have
$$
\lim_{i \to \infty}\H^{n-1}\bigl((\{\theta_i=a\}\cup \{\theta_i=b\})\cap B_{2R}(0)\bigr)=0
$$
because, by the coarea formula,
\begin{multline*}
\int_{-\infty}^{\infty}\H^{n-1}(\{\theta_i=s\}\cap B_{2R}(0))ds  =\int_{L^i\cap B_{2R}(0)}|H|\a\\
 \leq \sqrt{D_0R^n}\left(\int_{L^i\cap B_{2R}(0)}|H|^2\a\right)^{1/2}.
\end{multline*}
The first variation formula yields, for any vector field $Y$  supported in $B_{2R}(0)$
$$
\delta N^{i}(Y)=-\int_{N^{i}\cap B_{2R}(0)}\langle H,Y\rangle\a +\oint_{\partial N^i\cap B_{2R}(0)}\langle
Y,\nu\rangle\b,
$$
where $\nu$ denotes the exterior unit normal. Hence, whenever the sup norm of Y satisfies
$|Y|_{\infty}\leq 1$, we get
\begin{multline*}
|\delta N^{i}(Y)|\leq \sqrt{C_0R^n}\left(\int_{N^{i}\cap B_{2R}(0)}|H|^2\a\right)^{1/2}\\
+\H^{n-1}((\{\theta_i=a\}\cup \{\theta_i=b\})\cap B_{2R}(0))
\\+\H^{n-1}(\partial L^i\cap B_{2R}(0))
\end{multline*}
 Furthermore, if $\vartheta$ is any $n-1$
form compactly supported in $B_{2R}(0)$ with $|\vartheta|\leq 1$, then
\begin{align*}
|\partial N^i(\vartheta)|\leq & \H^{n-1}((\{\theta_i=a\}\cup \{\theta_i=b\})\cap B_{2R}(0))\\
&+\H^{n-1}(\partial L^i \cap B_{2R}(0)).
\end{align*}
We can now apply Allard compactness theorem for varifolds and Federer and Fleming compactness theorem for
currents (see \cite[Theorem 27.3]{Leon}) in order to complete the proof of the lemma.
\end{proof}

Condition (a) implies the existence of a finite set $F \subset\N$  such that, whenever $l \notin F$, we have for
all $i \in \N$
$$\{|\theta_i-(\bar\theta_1+l\pi)|\leq \pi\}\cap B_{2R}(0)=\emptyset.$$
\begin{lemm}
There is a universal constant $\gamma_n$ so that, for  all $\varepsilon<\pi/2,$
$$
 \lim_{i \to \infty}\sum_{l \in F}\H^n(\{|\theta_i-(\bar\theta_1+l\pi)|\leq \varepsilon\}\cap
   B_{\rho}(p))=\H^n (L \cap B_{\rho}(p))\geq \gamma_n \rho^n.
$$
\end{lemm}
\begin{proof}
The first equality is true because for almost all intervals $[a,b]$ such that
$$[a,b]\cap \{\bar\theta_1+l\pi\,\vert\,l\in \Z\}=\emptyset,$$
we have
$$
\limsup_{i \to \infty}\H^{n}\bigl(\{a \leq\theta_i\leq b\}\cap B_{\rho}(p)\bigr)=0.
$$
Otherwise we could, by Lemma \ref{lemma0}, extract a subsequence converging to a integer rectifiable varifold
$N$ with support in $L$ and such that $$\mu(B_{\rho}(p))>0,$$ where $\mu$ is the Radon measure associated to
$N$. This is impossible because for some positive $\delta$ we have
$$\sup_{\{a \leq\theta_i\leq b\}}|\cos(\theta_i-\bar\theta_1)|\leq 1-\delta,
$$
and so varifold convergence implies that
$$(1-\delta)\mu(B_{\rho}(p))\geq\lim_{i \to \infty}\int_{\{a \leq\theta_i\leq b\}\cap B_{\rho}(p)}
\Big|\re\left(e^{-i\bar\theta_1}\Omega\right)\Big|\a=\mu(B_{\rho}(p)).$$
\end{proof}

Renaming $\bar\theta_1$ to be $\bar\theta_1+l\pi$ for some $l$ in $F$, we can find a sequence $(\varepsilon_k)$
converging to zero and a constant $K=K(D_0)$ such that
\begin{equation}\label{lower}
\limsup_{i \to \infty}\H^n\bigl(\{|\theta_i-\bar\theta_1|\leq \varepsilon_k\}\cap B_{\rho}(p) \bigr)\geq K\rho^n
\end{equation}
for all $k \in \N$.

Applying Lemma \ref{lemma0} to
$$N^{i,k}\equiv\{|\theta_i-\bar\theta_1|\leq \varepsilon_k\},$$
we obtain two sequences $(N^k)$ and $(\widehat{N}^k)$ of stationary integer rectifiable varifolds and integral
currents with no boundary respectively. Its Radon measures are denoted by $\mu_k$ and $\widehat{\mu}_k$
respectively. Federer and Fleming compactness Theorem implies that $(\widehat{N}^k)$ has a subsequence that
converges in $B_{2R}(0)$ to an integral Lagrangian current ${L}_1$ with no boundary. Moreover,  $L_1$ is an
integral Special Lagrangian because it is calibrated by
$$
\vartheta\equiv\re \left(e^{-i\bar\theta_1}\Omega\right)
$$
and it is nonempty because, using \eqref{lower}, we obtain that for every nonnegative smooth $\phi$ compactly
supported in $B_{2R}(0)$
\begin{multline*}
\int_{L_1}\phi \a \geq \lim_{k\to\infty}\widehat{\mu}_k(\vartheta
\phi)=\lim_{k\to\infty}\lim_{i\to\infty}\int_{N^{i,k}}\vartheta\phi\\
 \geq \lim_{k\to\infty}\lim_{i\to\infty}\int_{N^{i,k}}\cos\varepsilon_k \phi\a =
\lim_{k\to\infty}\cos\varepsilon_k \mu_k(\phi)\geq K\rho^n.
\end{multline*}
Furthermore, the support of each integral current $\widehat{N}^k$ is a stationary rectifiable varifold which,
combined with the fact that
$$\widehat\mu_{k+1}(\phi)\leq \widehat\mu_{k}(\phi)$$
for every nonnegative $\phi$ compactly supported in $B_{2R}(0)$ and every $k \in \N$, implies that, for all $k$
sufficiently large, $\widehat{N}^k$ must coincide with ${L}_1$ in $B_{R}(0)$.
\end{proof}

Before proving Theorem A, we recall the monotonicity formula, found by Huisken in \cite{huisken}, valid for any
smooth family of $k$-dimensional submanifolds $(N_t)_{t\geq 0}$ moving by mean curvature flow in $\R^m$.
Consider the backward heat kernel
$$\Phi_{x_0,T}(x,t)=\frac{1}{(4\pi(T-t))^{k/2}}e^{-\frac{|x-x_0|^2}{4(T-t)}}.$$
When $(x_0,T)=(0,0)$, we denote it simply by $\Phi$. The following formula holds
$$\frac{d}{dt}\int_{N_t} f_t \Phi_{x_0,T}\a = \int_{N_t}
\left(\frac{d}{dt} f_t -\Delta
  f_t-\left|H+\frac{(\x-\x_0)^{\bot}}{2(T-t)}\right|^2 f_t\right)\Phi_{x_0,T}\a,$$
where $f_t$ is a smooth function with polynomial growth at infinity and $(\x-\x_0)^{\bot}$ denotes the
orthogonal projection on $(T_x N)^{\bot}$ of the vector determined by the point $(x-x_0)$ in $\R^m$.

Let $(L_t)_{0\leq t<T}$ be a solution to Lagrangian mean curvature flow  with a singularity at time $T$.
\begin{thma}
If $L_0$ is zero-Maslov class with bounded Lagrangian angle, then for any sequence of rescaled flows
$(L^i_s)_{s<0}$ at a singularity, there exist a finite set $\{\bar\theta_1,\ldots,\bar\theta_N\}$ and integral
Special Lagrangian cones
$$L_1,\ldots,L_N$$
such that, after passing to a subsequence, we have for every smooth function $\phi$ compactly supported, every
$f$ in $C^2(\R)$, and every $s<0$
$$
\lim_{i \to \infty}\int_{L^i_s}f(\theta_{i,s})\phi\a=\sum_{j=1}^N m_j f(\bar\theta_j)\mu_j(\phi),
$$
where $\mu_j$ and $m_j$ denote the Radon measure of the support of $L_j$ and its multiplicity respectively.

Furthermore, the set $\{\bar\theta_1,\ldots, \bar\theta_N\}$ does not depend on the sequence of rescalings
chosen.
\end{thma}

\begin{proof}
We start with the following lemma
\begin{lemm}\label{soul}
For any $a<b<0$ and any $R>0$, we have
$$ \lim_{i \to \infty}\int_a^b\int_{L^i_{s}\cap B_R(0)}\left(\big|\x^{\bot}\big|^2+|H|^2\right)\a ds=0.$$
\end{lemm}
\begin{proof}
From Huisken's monotonicity formula we have that, for all $i\in\N$,
$$
 \frac{d}{ds}\int_{L^i_s} \theta_{i,s}^2 \Phi\a = \int_{L_{s}^i}
 \left(-2|H|^2-\left|H-\frac{\x^{\bot}}{2s}\right|^2 \theta_{i,s}^2\right)\Phi\a
$$
and
$$
  \frac{d}{ds}\int_{L^i_s} \Phi\a =\int_{L_{s}^i} -\left|H-\frac{\x^{\bot}}{2s}\right|^2 \Phi\a.
$$
Using the scale invariance properties of the backward heat kernel, we obtain that
$$
 \lim_{i\to\infty}2\int_{a}^{b}\int_{L^i_s}|H|^2\Phi\a ds \leq\lim_{i \to \infty}
\left(\int_{L^i_{a}}\theta^2_{i,a}\Phi\a-\int_{L^i_{b}}\theta_{i,b}^2\Phi\a\right)=0 $$ and
$$
\lim_{i \to \infty}\int_{a}^{b}\int_{L^i_s}\left|H-\frac{\x^{\bot}}{2s}\right|^2\Phi\a ds=\lim_{i \to
\infty}\left(\int_{L^i_{a}}\Phi\a-\int_{L^i_{b}}\Phi\a\right)=0.
$$
Therefore
\begin{multline*}
\lim_{i \to \infty}\int_{a}^{b}\int_{L^i_s}\left|\frac{\x^{\bot}}{2s}\right|^2\Phi\a ds\\
\leq\lim_{i \to \infty} \int_{a}^{b}\int_{L^i_s}\left(\left|H-\frac{\x^{\bot}}{2s}\right|^2+|H|^2\right)\Phi\a
ds=0
\end{multline*}
and so the result follows.
\end{proof}
Pick $a<0$ for which
$$
\lim_{i \to \infty}\int_{L^i_{a}\cap B_R(0)}\left(\big|\x^{\bot}\big|^2+|H|^2\right)=0
$$
for all positive $R$.

The maximum principle implies that the Lagrangian angle $\theta_t$ is uniformly bounded and hence, by scale
invariance, the same is true for the Lagrangian angle of $L^i_{a}$.  Lemma \ref{bound} implies the existence of
a constant $D_0$ for which
$$
\H^n\left(L^i_{a}\cap B_R(0)\right)\leq D_0 R^n
$$
for all positive $R$. We can, therefore, apply Proposition \ref{general} to the sequence $(L^i_{a})$ and, after
a simple diagonalization argument, obtain a subsequence for which there are integral Special Lagrangian currents
$$L_1,\ldots,L_N$$
and a finite set $\{\bar\theta_1,\ldots,\bar\theta_N\}$ such that, for every smooth function $\phi$ compactly
supported  and every $f$ in $C^2(\R)$,
$$
\lim_{i \to \infty}\int_{L^i_{a}}f(\theta_{i,a})\phi\a=\sum_{j=1}^N m_j f(\bar\theta_j)\mu_j(\phi),
$$
where $\mu_j$ and $m_j$ denote the Radon measure and the multiplicity of $L_j$ respectively. The fact that
$$
\lim_{i \to \infty}\int_{L^i_{a}\cap B_R(0)}\big|\x^{\bot}\big|^2\a=0
$$
for all positive $R$ implies that the Special Lagrangians $L_j$ are all cones.

Next, we want to show that, for all $b<0$,
\begin{align*}
\lim_{i \to \infty}\int_{L^i_{b}}f(\theta_{i,b})\phi\a=\lim_{i \to
\infty}\int_{L^i_{a}}f(\theta_{i,a})\phi\a=\sum_{j=1}^N m_j f(\bar\theta_j)\mu_j(\phi).
\end{align*}
This comes from
\begin{multline*}
\frac{d}{ds}\int_{L^i_s}f(\theta_{i,s})\phi\a  = \int_{L^i_s}f^{\prime}(\theta_{i,s})\Delta \theta_{i,s}\phi\a
\\+\int_{L^i_s}f(\theta_{i,s})\langle H,D\phi\rangle\a -\int_{L^i_s}f(\theta_{i,s})|H|^2\phi\a
\end{multline*}
because, after integration with respect to the $s$ variable, all terms on the right hand side vanish when $i$
goes to infinity. We check this for the first term. Integrating by parts (and assuming $a<b$ for simplicity), we
obtain
\begin{multline*}
\int_{a}^{b}\int_{L^i_s}f^{\prime}(\theta_{i,s})\Delta \theta_{i,s}\phi\a
ds=-\int_{a}^{b}\int_{L^i_s}f^{\prime\prime}(\theta_{i,s})|\nabla \theta_{i,s}|^2\phi\a ds
\\
-\int_{a}^{b}\int_{L^i_s}f^{\prime}(\theta_{i,s})\langle\nabla \theta_{i,s},D\phi\rangle\a ds
\end{multline*}
and hence, by H\"olders's inequality, there is a constant $C=C(\phi,f,D_0,a,b)$ such that, for all $i\in \N$,
$$
\int_{a}^{b}\int_{L^i_s}\left|f^{\prime\prime}(\theta_{i,s})|\nabla \theta_{i,s}|^2\phi\right|\a ds \leq
C\int_{a}^{b}\int_{L^i_s}|H|^2\Phi\a ds
$$
and
$$
 \int_{a}^{b}\int_{L^i_s}\left|f^{\prime}(\theta_{i,s})\langle\nabla
\theta_{i,s},D\phi\rangle\right|\a ds \leq C\left(\int_{a}^{b}\int_{L^i_s}|H|^2\Phi\right)^{1/2}\a ds.
$$

Finally, we show that $\{\bar\theta_1,\ldots, \bar\theta_N\}$ does not depend on the sequence of rescalings
chosen. Let
$$\bigl(\widehat L^k_s\bigr)_{s<0}$$ be another sequence of rescaled flows for which there are Special Lagrangian
cones
$$\widehat L_1,\ldots, \widehat L_P$$
and a finite set $\bigl\{\hat\theta_1,\ldots,\hat\theta_P\bigr\}$ such that,  for every smooth function $\phi$
compactly supported, every $f$ in $C^2(\R)$, and every $s<0$
$$
\lim_{k \to \infty}\int_{\widehat L^k_s}f(\theta_{k,s})\phi\a=\sum_{j=1}^P \widehat m_j
f(\hat\theta_j)\widehat\mu_j(\phi),
$$
where $\widehat\mu_j$ and $\widehat m_j$ denote the Radon measure of the support of $L_j$ and its multiplicity
respectively.

For any real number $y$ and any integer $q$, we have the following evolution equation
$$\frac{d}{dt}(\theta_t-y)^{2q}=\Delta(\theta_t-y)^{2q}-2q(2q-1)(\theta_t-y)^{2p-2}|H|^2.$$
Applying the monotonicity formula to $(\theta_t-y)^{2q}$, we get that
$$\frac{d}{dt}\int_{L_t}(\theta_t-y)^{2q}\Phi_{x_0,T}\a\leq 0$$
and thus, by scale invariance, we obtain  for any  $s, \bar s<0$
\begin{multline*}
\lim_{i \to\infty}\int_{L^i_{s}}(\theta_{i,s}-y)^{2q}\Phi\a = \lim_{k
  \to\infty}\int_{\widehat L^k_{\bar s}}(\theta_{i,\bar s}-y)^{2q}\Phi\a\\
  =\lim_{t \to T}\int_{L_t}(\theta_s-y)^{2q}\Phi_{x_0, T}\a.
\end{multline*}
Therefore
$$\sum_{j=1}^N m_j (\bar\theta_j-y)^{2q}\mu_j(\Phi)=\sum_{j=1}^P \widehat m_j \bigl(\hat\theta_j-y\bigr)^{2q}
\widehat\mu_j(\Phi)$$ for all positive integer $q$ and all $y$ in $\R$ and this implies that
$$
\bigl\{\theta_1,\ldots\theta_N\bigr\}=\bigl\{\hat\theta_1,\ldots\hat\theta_P\bigr\}.
$$

\end{proof}

\section{Evolution Equations}\label{basic}

Let $L_0$ be a rational and zero-Maslov Lagrangian submanifold of $\C^n$. We will argue now that the rational
condition is preserved by the flow. Denoting by $F_t$ the normal deformation by mean curvature, we have
\begin{align*}
\frac{d}{dt}\int_{F_t(\gamma)}\lambda & = \frac{d}{dt}\int_{\gamma}F_t^*\lambda=\int_{\gamma}\L_H F_t^*\lambda\\
 & =\int_{\gamma}dF_t^*(H\lrcorner\lambda)+F_t^*(H\lrcorner 2\omega)=
\int_{\gamma}dF_t^*(H\lrcorner\lambda-2\theta_t)=0
\end{align*}
for every $[\gamma]$ in $H_1(L_0)$. Hence $$[\lambda]=[F_t^*(\lambda)] \quad\mbox{in }H^1(L_0)$$ for all times
where the solution exists smoothly and therefore it follows that
$$\lambda\left(H_1(L_t, \Z)\right)=\lambda\left(H_1(L_0, \Z)\right)= \{a2k\pi\,|\,k\in \Z\}.$$
Thus, there is a smooth family of multivalued functions
$$
\beta_{t}\,:\,L_t\longrightarrow \R/2\pi a \Z
$$
such that
$$
\nabla \beta_t(x)=(J\x)^{\top}\quad\mbox{for all }x \in L_t.
$$

\begin{prop}\label{evol}
The functions $\beta_t$ can be chosen so that
$$
\frac{d\beta_t}{dt}=\Delta\beta_t-2\theta_t.
$$
\end{prop}

\begin{proof}

Assume, without loss of generality, that the family of functions $\beta_t$ is smooth with respect to the time
parameter. We have
\begin{lemm}\label{evol1}
$$\Delta \beta_t=H\lrcorner\lambda,$$
\end{lemm}
\begin{proof}
 We use a normal coordinate system around the point $x$ and denote the coordinate vectors by
$\{\partial_1,\cdots,\partial_n\}$. The result follows from
\begin{align*}
\langle \nabla_{\partial_i} (J\x)^{\top},\partial_j \rangle & = \partial_i \langle J\x, \partial_j \rangle
-\langle (J\x)^{\top},D_{\partial_i}\partial_j\rangle = \langle J\partial_i, \partial_j \rangle + \langle
(J\x)^{\bot}, D_{\partial_i}\partial_j \rangle\\ & = \langle J\x, A_{ij} \rangle.
\end{align*}
\end{proof}
Thus,
\begin{align*}
d\left(\frac{d\beta_t}{dt}\right)=\frac{d \lambda}{dt} & = \L_H \lambda = d(H\lrcorner\lambda)+H\lrcorner
2\omega
  =  d (\Delta \beta_t-2\theta_t)
\end{align*}
and so we can add a time dependent constant to each $\beta_t$ so that the desired result follows.
\end{proof}

Given any $t_0$ in $\R$ and any $k$ in $\Z$, the function
$$
u_{t}\equiv\cos\left(\frac{k(\beta_t+2(t-t_0)\theta_t)}{a}\right)
$$
is well defined on $L_t$. If $L_0$ is exact, take $a=1$. A straightforward computation using Proposition
\ref{evol} and
$$J\x^{\bot}=(J\x)^{\top}$$
gives
\begin{cor}\label{beta}
$$
\frac{d u_{t}}{dt}=\Delta u_{t}+u_{t}\left|\frac{k(\x^{\bot}+2(t_0-t)H)}{a}\right|^2.
$$
\end{cor}

\section {Proof of Compactness Theorem B}\label{B}

\begin{thmb}
If $L_0$ is almost-calibrated  and rational, then, after passing to a subsequence of $(L^i_s)_{s<0}$, the
following property holds for all $R>0$ and almost all $s <0$. \par For any convergent subsequence (in the Radon
measure sense) $\Sigma^i$ of connected components of $B_{4R}(0)\cap L^i_s$ intersecting $B_R(0)$, there exists a
Special Lagrangian cone $L$ in $B_{2R}(0)$ with Lagrangian angle $\bar\theta$ such that
$$
\lim_{i \to \infty}\int_{\Sigma^i}f(\theta_{i,s})\phi\a=m f(\bar\theta)\mu(\phi),
$$
for every $f$ in $C(\R)$ and every smooth $\phi$ compactly supported in $B_{2R}(0)$, where $\mu$ and $m$ denote
the Radon measure of the support of $L$ and its multiplicity respectively.
\end{thmb}

\begin{proof}
The almost-calibrated condition is preserved by the flow and implies the following lemma.
\begin{lemm}\label{iso}
There is a constant $C_1$ such that, for  all $s<0$,
$$\left(\H^n(A)\right)^{(n-1)/n}\leq C_1 \H^{n-1}(\partial A),$$
where $A$ is any open subset of $L^i_s$ with rectifiable boundary.
\end{lemm}
% When $n=2$ the lemma will follow from the Michael-Simon Sobolev inequality
%(see \cite[Theorem 18.6]{Leon})
%$$
%\left(\H^2(A)\right)^{1/2}\leq C\int_A |H| +C\H^{1}(\partial A)
%$$
%for some universal constant $C$. In this case we have
%$$
%\left(\H^2(A)\right)^{1/2}\leq C\left(\H^2(A)\right)^{1/2} \left(\int_A |H|^2\right)^{1/2}+C\H^{1}(\partial A)
%$$
%and so we get the desired result whenever $$C^2 \int_{L^i_s\cap B_R(0)} |H|^2\leq 1/4.$$
\begin{proof}
The Isoperimetric Theorem \cite[Theorem 30.1]{Leon} guarantees the existence of an integral current $B$ with
compact support such that $\partial B=\partial A$ and for which
$$
\left(\H(B)\right)^{(n-1)/n}\leq C \H^{n-1}(\partial A),
$$
where $C=C(n)$. If $T$ denotes the cone over the current $A-B$ (see \cite[page 141]{Leon}), then $\partial
T=A-B$ and thus, because $$\re \,\Omega_{|L^i_s}=\cos \theta_{i,s} \geq \varepsilon_0$$ for some positive
$\varepsilon_0$, we obtain
\begin{align*}
\H^{n}(A) & \leq \varepsilon_0^{-1}\int_A \re\, \Omega=\varepsilon_0^{-1}\int_B \re\, \Omega+\partial T(\re\, \Omega)\\
& \leq \varepsilon_0^{-1}\H^n(B)+T(d\re\,\Omega)\leq \varepsilon_0^{-1}\left(C\H^{n-1}(\partial
A)\right)^{n/(n-1)}.
\end{align*}
\end{proof}

The discussion in Section \ref{basic} implies the existence of $a\in\R$ and of a family of multivalued functions
$$\beta_{i,s}\,:\,L^i_s\longrightarrow\R/\sigma_i^2a2\pi\Z$$% k\,|\,k\in\Z\}$$
such that
$$
\nabla \beta_{i,s}(x)=(J\x)^{\top}$$ for all $x \in L^i_s$ and all $s<0$. Furthermore, we can choose a bounded
sequence $(b_i)$ so that, for any real number $s_0,$
$$u_{i,s}\equiv\cos\left(\frac{\beta_{i,s}+2(s-s_0)\theta_{i,s}}{b_i}\right)$$
is a well defined function. After passing to a subsequence, the sequence $(b_i)$ converges to $b\neq 0$ and, for
simplicity, we assume that $b=1$. Furthermore, from Lemma \ref{soul}, we can also  assume that
$$
\lim_{i \to \infty}\int_{L^i_{-1}\cap
        B_{R}(0)}\left(|H|^2+\big|\x^{\bot}\big|^2\right)\a=0
$$
for all $R>0.$

\begin{lemm}\label{lemma1} There is a set
  $$
\{(\cos\bar\beta_1,\sin\bar\beta_1) ,\ldots,(\cos\bar\beta_Q,\sin\bar\beta_Q)\}
$$
 and integral Special Lagrangian cones
$$P_1,\ldots,P_Q$$ such that, after passing to a subsequence,
we have for all smooth $\phi$ with compact support and all $f$ in $C(\R)$,
\begin{align*}
 & \lim_{i\to\infty}\int_{L^i_{-1}}f(\cos(\beta_{i,-1}/b_i))\phi\a=\sum_{k=1}^Q
p_kf(\cos \bar\beta_k )\nu_{k}(\phi)\\
 & \lim_{i\to\infty}\int_{L^i_{-1}}f(\sin(\beta_{i,-1}/b_i))\phi\a=\sum_{k=1}^Q
p_kf(\sin \bar\beta_k)\nu_{k}(\phi),
\end{align*}
 where $\nu_k$ and the positive integer $p_k$ denote the Radon measure of the support of $P_k$ and
 its multiplicity respectively.
\end{lemm}

\begin{proof}
Let $(R_k)$ denote a sequence of positive numbers going to infinity. We start by arguing the existence of a
uniform bound on the number of connected components of $L^i_{-1}\cap B_{4R_k}(0)$ that intersect $B_{R_k}(0)$.
For any $x$ in $L^i_{-1}\cap B_{2R_k}(0)$, denote the intrinsic ball of radius $r$ around $x$ by
$\widehat{B}_{i}(x,r)$. Set
$$\psi_i(r)\equiv\H^n\left(\widehat{B}_{i}(x,r)\right)$$ which has, for almost all $r$, derivative given by
$$\psi_i^{\prime}(r)=\H^{n-1}\left(\partial \widehat{B}_{i}(x,r)\right).$$
We know from Lemma \ref{iso}  that, for all $r<R_k$,
$$(\psi_i(r))^{(n-1)/n}\leq C_1\psi_i^{\prime}(r)
$$ and so
$$\H^n\left(\widehat{B}_{i}(x,r)\right)\geq Kr^n$$ for all $x$ in
$L^i_{-1}\cap B_{2R_k}(0)$, where $K=K(n,C_1)$. Hence, each connected component has area bigger than $KR^n$ and
so the claim follows from the uniform area bounds for $L^i_{-1}$ (Lemma \ref{bound}).

From Proposition \ref{general} we know that, after passing to a subsequence, all the connected components of
$L^i_{-1}\cap B_{4R_k}(0)$ intersecting $B_{R_k}(0)$ converge to a union of Special Lagrangian cones in
$B_{2R_k}(0)$. Moreover,
$$|\nabla \beta_{i,-1}(x)|=\big|(J \x )^{\top}\big|=\big|\x^{\bot}\big|$$
and thus the functions $\cos(\beta_{i,-1}/b_i)$ and $\sin(\beta_{i,-1}/b_i)$ satisfy the conditions of
Proposition \ref{cont}. We can, therefore, apply this result to all the  connected components of $L^i_{-1}\cap
B_{2R_k}(0)$ intersecting $B_{R_k}(0)$. A standard diagonalization method finds a subsequence that works for all
$R_k$ and so the lemma is proved.
\end{proof}

Combining this lemma with Theorem A we obtain that, after a rearrangement of the supports of the Special
Lagrangian cones and its multiplicities (which we still denote by
$$L_1,\ldots,L_N$$ and $m_1,\ldots,m_N$ respectively),  we can
assume that for all $\phi$ with compact support, all $f$ in $C(\R)$, and all $y \in \R$,
 $$
 \lim_{i\to \infty}\int_{L^i_{-1}}f\left(\cos \left(\frac{\beta_{i,-1}+2y\theta_{i,-1}}{b_i}\right)\right)\phi\a=
 \sum_{j=1}^N m_j f(\cos(\bar\beta_j+2y\bar\theta_j))\mu_j(\phi)
 $$
where $\mu_j$ denotes the Radon measure of the support of $L_j$ and the elements of the set $$\{(\cos
\bar\beta_1,\sin\bar\beta_1,\bar\theta_1),\ldots,(\cos \bar\beta_N,\sin\bar\beta_Q,\bar\theta_N)\}$$ are all
distinct.

Using the evolution equation for $u_{i,s}$ we show
\begin{lemm}\label{lemma2}
For all $\phi$ with compact support, all $f$ in $C^2(\R)$, and all $ s <0$,
$$
\lim_{i\to \infty}\int_{L^i_s}f(\cos(\beta_{i,s}/b_i))\phi\a= \sum_{j=1}^N m_j
f(\cos(\bar\beta_j-2(s+1)\bar\theta_j))\mu_j(\phi).
$$
\end{lemm}

\begin{proof}
 Corollary \ref{evol} implies that, for all $\phi$ with compact support, all $f$ in $C^2(\R),$ and all $s_0<0,$
\begin{multline}\label{silence}
\frac{d}{ds}\int_{L^i_s}f(u_{i,s})\phi\a = \int_{L^i_s}f^{\prime}(u_{i,s})\Delta u_{i,s}\phi\a\\ +
\int_{L^i_s}f^{\prime}(u_{i,s})u_{i,s}\left|\frac{\x^{\bot}+2(s_0-s)H}{b_i}\right|^2\phi\a
+\int_{L^i_s}f(u_{i,s})\langle H,D\phi\rangle\a\\- \int_{L^i_s}f(u_{i,s})|H|^2\phi\a.
\end{multline}
From Lemma \ref{soul}, we obtain that (assuming $-1<s_0<0$ for simplicity)
\begin{multline*}
\lim_{i\to \infty}\int_{-1}^{s_0}\int_{L^i_s\cap B_R(0)}\left|\frac{\x^{\bot}+2(s_0-s)H}{b_i}\right|^2\a
\\
\leq\lim_{i\to\infty}8\int_{-1}^{s_0} \int_{L^i_s\cap
B_R(0)}\left(\frac{(s-s_0)^2|H|^2+\left|\x^{\bot}\right|^2}{b_i^2}\right)\a=0
\end{multline*}
for all positive $R$.

This inequality allows us to argue in the same way as it was done in the proof of Theorem A and show that, after
integration with respect to the $s$ variable, all terms on the right hand side of \eqref{silence} converge to
zero when $i$ goes to infinity. Thus, because
$$u_{i,s_0}=\cos(\beta_{i,s_0}/b_i)\quad\mbox{and}\quad u_{i,-1}=\cos \left(\frac{\beta_{i,-1}
-2(1+s_0)\theta_{i,-1}}{b_i}\right),
$$ we obtain from Lemma \ref{lemma2}
\begin{multline*}
\lim_{i\to \infty}\int_{L^i_{s_0}}f(\cos(\beta_{i,s_0}/b_i))\phi\a \\= \lim_{i\to
\infty}\int_{L^i_{-1}}f\left(\cos \left(\frac{\beta_{i,-1}-2(1+s_0)\theta_{i,-1}}{b_i}\right)\right)\phi\a\\
=\sum_{j=1}^N m_j f(\cos(\beta_j-2(1+s_0)\theta_j))\mu_j(\phi).
\end{multline*}
The result follows from the arbitrariness of $s_0$.
\end{proof}

The proof of the theorem can now be completed. Because the elements of the set $$\{(\cos
\bar\beta_1,\sin\bar\beta_1,\bar\theta_1),\ldots,(\cos \bar\beta_N,\sin\bar\beta_Q,\bar\theta_N)\}$$ are all
distinct, we get that, for all but countably many $s$, the real numbers
$$
\cos(\bar\beta_1-2(s+1)\bar\theta_1),\ldots,\cos(\bar\beta_N-2(s+1)\bar\theta_N)
$$are all distinct. Moreover, Lemma \ref{soul} implies that, for almost all $s<0$,
$$
\lim_{i \to \infty}\int_{L^i_{s}\cap
        B_{R}(0)}\left(|H|^2+\big|\x^{\bot}\big|^2\right)\a=0
$$
for all $R>0$.

Pick $s$ so that both conditions described above hold and consider a subsequence of connected components
$\Sigma^i$ of $B_{4R}(0)\cap L^i_s$ intersecting $B_R(0)$ that converges weakly to $\Sigma$. The arguments
presented in the proof of Lemma \ref{soul} imply that $\Sigma$ has positive measure. We first show that $\Sigma$
is a Special Lagrangian cone.

Proposition \ref{cont} can be applied to the sequence $\Sigma_i$ and thus, after passing to a subsequence,
$(\cos (\beta_{i,s}/b_i))$ converges  to a constant $\gamma$. Define $f\in C^2(\R)$ to be a nonnegative cutoff
function that is one in small neighborhood of $\gamma$ and zero everywhere else.

Denoting by $\mu_{\Sigma}$  the Radon measure of $\Sigma$, we obtain from Lemma \ref{lemma2} that for every
nonnegative test function $\phi$ with support in $B_{2R}(0)$
\begin{multline*}
\mu_{\Sigma}(\phi)=\lim_{i\to\infty}\int_{\Sigma_i}\phi\a=\lim_{i\to\infty}\int_{\Sigma_i}f(\cos (\beta_{i,s}/b_i))\phi\a \\
\leq\lim_{i\to\infty}\int_{L_i}f(\cos (\beta_{i,s}/b_i))\phi\a=\sum_{j=1}^N m_j
f(\cos(\bar\beta_j-2(s+1)\bar\theta_j))\mu_j(\phi).
\end{multline*}
Because the support of $f$ can be chosen arbitrarily small and the real numbers
$$\cos(\bar\beta_1-2(s+1)\bar\theta_1),\ldots,\cos(\bar\beta_N-2(s+1)\bar\theta_N)$$
are all distinct, the above inequality implies that
$$\gamma=\cos(\bar\beta_{j_0}-2(s+1)\bar\theta_{j_0})$$ for a unique $j_0$. Thus
$$
\mu_{\Sigma}(\phi)\leq m_{j_0} \mu_{j_0}(\phi)
$$
for every $\phi\geq 0$ and, as a result, the support of $\Sigma$ must be contained in $L_{j_0}$.

 Finally, suppose there are $f$ continuous and $\phi$ compactly supported in $B_{2R}(0)$ such that
$$
  \int_{\Sigma^i}f(\theta_{i,s})\phi\a
$$
has two distinct convergent subsequences. We can use  Proposition \ref{general} to get a contradiction because
$L_0$ being almost-calibrated implies that any two Special Lagrangian cones with support contained in the
support of $\Sigma$ have the same Lagrangian angle.
\end{proof}

\appendix

\section{ }
Suppose we have a sequence of functions $(\alpha_i)$ defined on a  sequence of manifolds $(N^i)$ converging
weakly to $N$ and such that the $L^2$-norm of $|\nabla \alpha_i|$ converges to zero. The next proposition gives
conditions under which, after passing to a subsequence, $(\alpha_i)$ converges to a constant. Before giving its
proof, we comment on the necessity of all the hypothesis.
\begin{prop}\label{cont}
Let $(N^i)$ and $(\alpha_i)$ be a sequence of smooth $k$-submanifolds in $\R^n$  and  smooth functions on $N^i$
respectively, such that $(N^i)$ converges weakly to an integer rectifiable stationary $k$-varifold $N$. We
assume that, for some $R>0$, the following properties hold:
\begin{itemize}
\item[a)] There exists a constant $D_0$ such that  $$\H^k(N^i\cap B_{3R}))\leq D_0R^k$$ for all $i\in\N$, and
$$\left(\H^k(A)\right)^{(k-1)/k}\leq D_0 \H^{k-1}(\partial A)$$
for all open subsets $A$ of $N^i\cap B_{3R}$ with rectifiable boundary.

\item[b)]$$\lim_{i \to \infty}\int_{N^i\cap B_{3R}(0)}\left(|H|^2+|\nabla \alpha_i|^2\right)\a=0.$$

\item[c)] There exists a constant $D_1$ for which
$$\sup_{N^i\cap B_{3R}(0)}|\nabla \alpha_i|+R^{-1}\sup_{N^i\cap B_{3}(0)}|\alpha_i|\leq D_1$$ for all $i \in \N$.

\item[d)] For all $i \in\N$, $$N^i\cap B_{2R}(0)\quad\mbox{ is connected}$$ and
$$ \partial(N^i\cap B_{3R}(0))\subset\partial B_{3R}(0).$$
\end{itemize}
Then, there is  a real number $\alpha$ such that, after passing to a subsequence, we have for all $\phi$ with
compact support in $B_R(0)$ and all $f$ in $C(\R)$
 $$\lim_{i \to
   \infty}\int_{N^i}f(\alpha_i)\phi=f(\alpha)\mu_N(\phi),$$
where $\mu_N$ denotes the Radon measure associated to $N$.
\end{prop}

 The first hypothesis is needed in order to ensure lower density bounds on $N^i$. The third hypothesis is
 essential because, without the pointwise bounds on $|\nabla \alpha_i|$ and $\alpha_i$, the result would be false.
 Finally, the last hypothesis is needed because otherwise the proposition would fail for trivial reasons.

\begin{proof} It suffices to
find $\alpha \in \R$ and a sequence $(\varepsilon_j)$ converging to zero such that, for some appropriate
subsequence, we have for all $j \in \N$
$$\lim_{i \to \infty}\H^k(\{|\alpha_i-\alpha|\leq \varepsilon_j\}\cap B_{R}(0))=\H^k(N\cap B_R(0)).$$

For the rest of this proof, $K=K(D_0,D_1,k)$ will denote a generic constant depending only on the mentioned
quantities. Choose any sequence $(x_i)$ in $N^i\cap B_R(0)$. After passing to a subsequence, we have that $$
\lim_{i \to\infty}x_i=x_0\quad\mbox{and}\quad\lim_{i \to
  \infty}\alpha_i(x_i)=\alpha
$$ for some $x_0 \in B_R(0)$ and $\alpha \in \R$. Furthermore, consider also
a sequence $(\varepsilon_j)$ converging to zero such that, for all $j\in \N$,
$$
\lim_{i \to \infty}\H^{k-1}\bigl( \{\alpha_i=\alpha\pm\varepsilon_j\}\cap B_{3R}\bigr)=0.
$$
Such a subsequence exists because, by the coarea formula, we have
\begin{multline*}
\lim_{i \to \infty}\int_{-\infty}^{\infty}\H^{k-1}\bigl( \{\alpha_i=s\}\cap B_{3R}\bigr)ds
=\lim_{i \to \infty}\int_{N^i\cap B_{3R}}|\nabla \alpha_i|\a\\
\leq \lim_{i \to \infty}K R^{k/2}\left(\int_{N^i\cap B_{3R}}|\nabla \alpha_i|^2\a\right)^{1/2} =0.
\end{multline*}
Define $$
 N^{i,\alpha,j}\equiv\{|\alpha_i-\alpha|\leq \varepsilon_j\}.
$$
The first variation formula yields, for any vector field $Y$ supported in $B_{3R}$,
$$
\delta N^{i,\alpha,j}(Y)=-\int_{N^{i,\alpha,j}\cap B_{2R}}\langle H,Y\rangle\a
+\oint_{\partial\{|\alpha_i-\alpha|\leq \varepsilon_j\}\cap
  B_{2R}}\langle Y,\nu\rangle\b
$$
where $\nu$ denotes the exterior unit normal. Hence, whenever the sup norm of Y satisfies $|Y|_{\infty}\leq 1$,
we get
\begin{multline*}
|\delta N^{i,\alpha,j}(Y)|\leq KR^{k/2}\left(\int_{N^{i,\alpha,j}\cap B_{2R}}|H|^2\a
\right)^{1/2}\\+\H^{k-1}\bigl(\{\alpha_i=\alpha\pm\varepsilon_j\}\cap B_{2R}\bigr).
\end{multline*}
We can now apply Allard compactness theorem to conclude that, after passing to a subsequence, we have
convergence to an integer rectifiable stationary varifold $N^{\alpha,j}$. By a standard diagonalization
argument, we can find a subsequence that works for every positive integer $j$.

\begin{lemm}
For all $j \in \N$, $$ \H^k\bigl(N^{\alpha,j}\cap B_R(x_0)\bigr)\geq K R^k.
$$
\end{lemm}
\begin{proof}
Set
$$
\psi_i(s)\equiv\H^k\bigl(\{|\alpha_i-\alpha_i(x_i)|\leq s\}\cap B_s(x_i)\bigr)
$$
which, by the coarea formula, has derivative equal to
\begin{multline*}
\psi_i^{\prime}(s)=\oint_{\partial B_s(x_i)\cap \{|\alpha_i-\alpha_i(x_i)|\leq
  s\}}\frac{|x-x_i|}{|(\x-\x_i)^{\top}|}\b \\+\oint_{B_s(x_i)\cap
\partial \{|\alpha_i-\alpha_i(x_i)|\leq s\}}\frac{1}{|\nabla\alpha_i|}\b
\end{multline*}
for almost all $s$. We can estimate
\begin{multline*}
\psi_i^{\prime}(s)  \geq \H^{k-1}\bigl(\partial B_s(x_i)\cap \{|\alpha_i-\alpha_i(x_i)|\leq
 s\}\bigr)\\
  +K \H^{k-1}\bigl(B_s(x_i)\cap \partial \{|\alpha_i-\alpha_i(x_i)|\leq s\}\bigr)\\
  \geq K\H^{k-1}\bigl(\partial (B_s(x_i)\cap \{|\alpha_i-\alpha_i(x_i)|\leq s\})\bigr)
\end{multline*}
and so, using the isoperimetric condition a), we obtain
\begin{equation*}
 (\psi_i(s))^{(k-1)/k}\leq D_0 \H^{k-1}\bigl(\partial (B_s(x_i)\cap
\{|\alpha_i-\alpha_i(x_i)|\leq s\})\bigr)\leq K\psi_i^{\prime}(s)
\end{equation*}
for almost all $s\leq R$. This implies that
$$
s^{-k}\H^k\bigl(\{|\alpha_i-\alpha_i(x_i)|\leq s\}\cap B_s(x_i)\bigr)\geq K$$ for all $s\leq R$. This inequality
and the inclusion $$ \{|\alpha_i-\alpha_i(x_i)|\leq \varepsilon_j/2\}\cap
B_{\varepsilon_j/2}(x_i)\subset\{|\alpha_i-\alpha|\leq \varepsilon_j\}\cap B_{\varepsilon_j}(x_0),
$$ valid for all $i$ sufficiently large, imply that
\begin{multline*}
 \varepsilon_j^{-k}\H^k\bigl(N^{i,\alpha,j}\cap
 B_{\varepsilon_j}(x_0)\bigr)\\
 \geq \varepsilon_j^{-k}\H^k\bigl(\{|\alpha_i-\alpha(x_i)|\leq \varepsilon_j/2\}\cap
 B_{\varepsilon_j/2}(x_i)\bigr)\geq K
\end{multline*}
for all $i$ sufficiently large.
Taking the limit when $i$ goes to infinity and recalling that $N^{\alpha,j}$ is a stationary varifold we get, by
the monotonicity formula, that $$R^{-k}\H^k\bigl(N^{\alpha,j}\cap B_R(x_0)\bigl)\geq
\varepsilon_j^{-k}\H^k\bigl(N^{\alpha,j}\cap B_{\varepsilon_j}(x_0)\bigl)\geq K$$ for all $j\in \N$.
\end{proof}
Suppose that for some positive integer $j$ we have
$$\H^k\bigl(N^{\alpha,j}\cap B_R(0)\bigr)<\H^k(N\cap B_R(0)).$$
Repeating the same type of arguments, we can find $y_0$ in $B_R(0)$ and a closed interval $I$ disjoint from
$[\alpha-\varepsilon_j,\alpha+\varepsilon_j]$ so that, after passing to a subsequence,
$$
\lim_{i \to \infty}\H^k\bigl(\alpha_i^{-1}(I)\cap B_R(y_0)\bigr)\geq K R^k.
$$
Given any positive integer $p$, pick disjoint closed intervals
$$I_1,\cdots,I_p$$ lying between $I$ and
$[\alpha-\varepsilon_j,\alpha+\varepsilon_j]$. The connectedness of $N^i\cap B_{2R}(0)$ implies that all
$\alpha_i^{-1}(I_l)\cap B_{2R}(0)$ are nonempty for $i$ sufficiently large. Hence, arguing as before, we find
$y_1,\dotsc, y_p$ in $B_{2R}(0)$
 such that, after passing to a subsequence, $$
 \lim_{i \to \infty}\H^k\bigl(\alpha_i^{-1}(I_l)\cap B_R(y_l)\bigr)\geq K R^k,
$$
for all $l$ in $\{1,\dotsc,p\}$. This implies that
\begin{align*}
\lim_{i \to \infty}\H^k\bigl(N^i\cap B_{2R}(0)\bigr) & \geq \lim_{i \to
  \infty}\sum_{l=1}^p\H^k\bigl(\alpha_i^{-1}(I_j)\cap B_R(y_l)\bigr)\\
&\geq pKR^k.
\end{align*}
Choosing $p$ sufficiently large we get a contradiction.
\end{proof}

\section{}

The next lemma is a simple modification of a result that can be found in Ecker's book \cite{eckernotes} and
Ilmanen's preprint \cite{ilmanen}. The proof is the same but we write it here for the sake of completeness.

\begin{lemm}\label{bound}
Let $(M_t)_{t\ge 0}$ be family of $k$-dimensional submanifolds $(M_t)_{t\geq 0}$ moving by mean curvature flow
in $\R^m$. Assume there are constants $A_0$ and $R_0$ such that
$$ \H^k\bigl(M_0 \cap B_r(0)\bigr)\leq A_0 r^k,$$
for all $r \geq R_0$. Then, for all $t \geq t_0$ and $x_0 \in \R^m$, there is a constant
$C=C(A_0,R_0/\sqrt{t_0}, |x_0|)$ such that
$$
\H^k\bigl(M_t\cap B_r(x_0)\bigr)\leq C r^k
$$
for all $r>0$.
\end{lemm}
\begin{proof} In what follows, $C=C(A_0,t_0^{-1},R_0, |x_0|)$ will denote a constant depending only on the
mentioned quantities. Using the monotonicity formula we obtain
\begin{align*}
\frac{\H^k\bigl(M_t\cap B_r(x_0)\bigr)}{r^k} & \leq C\int_{M_t}\frac{1}{(4\pi
r^2)^{k/2}}e^{-\frac{|x-x_0|^2}{4r^2}}\a\\
& \leq C\int_{M_0}\frac{1}{(4\pi (t+r^2))^{k/2}}e^{-\frac{|x-x_0|^2}{4(t+r^2)}}\a\\
&\leq C\int_{M_0}\frac{1}{(4\pi
  (t+r^2))^{k/2}}e^{-\frac{|x|^2}{8(t+r^2)}}\a \\&\leq C\int_{\lambda M_0}e^{-|x|^2}\a,
\end{align*}
where $\lambda\equiv(8(t+r^2))^{-1/2}$. For all $s\geq \lambda R_0$ we have
$$
\H^k\bigl(\lambda M_0 \cap B_s(0)\bigr)\leq A_0 s^k
$$
and thus, setting $R_1\equiv\max\{2, (8t_0)^{-1/2}R_0\}$, the result follows from
\begin{multline*}
\int_{\lambda M_0}e^{-|x|^2}\a\leq A_0 R_1^k+\int_{\lambda
  M_0\setminus B_{R_1}}e^{-|x|^2}\a\\ =  A_0 R_1^k+\sum_{j\geq
  0}\int_{\lambda M_0\cap \left(B_{R_1^{j+1}}\setminus B_{R^j_1}\right)}e^{-|x|^2} \a\\
   \leq A_0 R_1^k+\sum_{j\geq 0}A_0 R_1^{j+1} e^{-R_1^{2j}}.
\end{multline*}
\end{proof}

\bibliographystyle{amsbook}

\begin{thebibliography}{99}
\bibitem {anciaux}
H. Anciaux, Mean curvature flow and self-similar submanifolds, {\bf S\'eminaire de Th\'eorie Spectrale et
G\'emom\'etrie. Vol. 21} Ann\'ee 2002--2003, 43--53.
\bibitem {ang} S. Angenent, Parabolic equations for curves on
      surfaces. II. Intersections,  blow-up and generalized solutions,
{\bf Ann. of Math. (2) 133} (1991), 171--215.
\bibitem {CL} J. Chen and J. Li, Singularity of mean curvature flow of Lagrangian submanifolds,
    {\bf Invent. Math. 156} (2004), 25--51.
\bibitem {CLT}  J. Chen, J. Li and G. Tian,
     Two-dimensional graphs moving by mean curvature flow,
     {\bf Acta Math. Sin. 18}, (2002), 209--224.
\bibitem {eckernotes} K. Ecker, Regularity theory for mean curvature flow,
{\bf Progress in Nonlinear Differential Equations and their Applications, 57}, Birkhäuser Boston, MA, 2004.

\bibitem {Ecker-Huisken} K. Ecker and G. Huisken,
Mean curvature evolution of entire graphs, {\bf Ann. of Math. (2) 130} (1989), 453--471.

\bibitem {HL} R. Harvey and H. B. Lawson, H.
Calibrated geometries, {\bf Acta Math. 148} (1982), 47--157.

\bibitem {huisken} G. Huisken,
Asymptotic behavior for singularities of the mean curvature flow, {\bf J. Differential Geom. 31} (1990),
285--299.
\bibitem {ilmanen} T. Ilmanen, Singularities of Mean Curvature Flow of Surfaces, preprint.
\bibitem {Leon} L. Simon,
Lectures on geometric measure theory, {\bf Proceedings of the Centre for Mathematical Analysis, Australian
National University, 3}.
\bibitem {schoen}R. Schoen and J.  Wolfson,
Minimizing area among Lagrangian surfaces: the mapping problem, {\bf J. Differential Geom. 58} (2001), 1--86.
\bibitem {smo0} K. Smoczyk,
A canonical way to deform a Lagrangian submanifold, preprint.
\bibitem {smo} K.~Smoczyk,
    Harnack inequality for the Lagrangian mean curvature flow,
    {\bf Calc. Var. Partial Differential Equations 8} (1999),
    247--258.
\bibitem{smo1}K. Smoczyk,
 Angle theorems for the Lagrangian mean curvature flow,
  {\bf Math. Z. 240 }(2002), 849--883.

\bibitem {smo2} K. Smoczyk,
    Longtime existence of the Lagrangian mean curvature flow,
   {\bf Calc. Var. Partial Differential Equations 20} (2004), 25--46.
\bibitem {SM} K. Smoczyk and M.-T. Wang,
     Mean curvature flows of Lagrangians submanifolds with convex
     potentials,
     {\bf J. Differential Geom. 62} (2002), 243--257.
\bibitem {Wa4} M.-P. Tsui and M.-T. Wang,
     Mean curvature flows and isotopy of maps between spheres,
     {\bf Comm. Pure Appl. Math. 57} (2004), 1110--1126.
\bibitem {Wa1} M.-T. Wang,
     Mean curvature flow of surfaces in Einstein four-manifolds,
     {\bf J. Differential Geom. 57} (2001), 301--338.

\bibitem {Wa0} M.-T. Wang,
Deforming area preserving diffeomorphism of surfaces by mean curvature flow, {\bf Math. Res. Lett. 8} (2001),
651--661.


\bibitem {Wa2} M.-T. Wang,
     Long-time existence and convergence of graphic mean curvature
     flow in arbitrary codimension,
    {\bf Invent. Math. 148} (2002), 525--543.

\bibitem {Wa3} M.-T. Wang,
Gauss maps of the mean curvature flow,   {\bf Math. Res. Lett. 10} (2003), 287--299.
   % \bibitem{Wa5} M.-T. Wang, A convergence result of the Lagrangian mean curvature flow, preprint.
\bibitem {white} B. White, A local regularity theorem for mean curvature flow, preprint.

\end{thebibliography}

\vspace{20mm}

\end{document}